\documentclass[11pt]{amsart} 
\usepackage{amsmath} 
\usepackage{amsxtra}
\usepackage{amscd} 
\usepackage{amsthm} 
\usepackage{amsfonts}
\usepackage{amssymb} 
\usepackage{eucal}

\newtheorem{cor}[subsection]{Corollary}
\newtheorem{lem}[subsection]{Lemma}
\newtheorem{prop}[subsection]{Proposition}

\newtheorem{thm}[subsection]{Theorem}

\theoremstyle{definition}

\theoremstyle{remark}

\newcommand{\thmref}[1]{Theorem~\ref{#1}}
\newcommand{\secref}[1]{Sect.~\ref{#1}}
\newcommand{\lemref}[1]{Lemma~\ref{#1}}
\newcommand{\propref}[1]{Proposition~\ref{#1}}
\newcommand{\corref}[1]{Corollary~\ref{#1}}

\newcommand{\nc}{\newcommand}
\nc{\renc}{\renewcommand}
\nc{\ssec}{\subsection}
\nc{\sssec}{\subsubsection}
\nc{\on}{\operatorname}

\nc\ol{\overline}
\nc\wt{\widetilde}
\nc\tboxtimes{\wt{\boxtimes}}
\nc{\alp}{\alpha}

\nc{\ZZ}{{\mathbb Z}}
\nc{\NN}{{\mathbb N}}
\nc{\CC}{{\mathbb C}}
\nc{\OO}{{\mathbb O}}
\renc{\SS}{{\mathbb S}}
\nc{\DD}{{\mathbb D}}
\nc{\GG}{{\mathbb G}}

\nc{\Fq}{{\mathbb F}_q}
\nc{\Fqb}{\ol{{\mathbb F}_q}}
\nc{\Ql}{\ol{{\mathbb Q}_\ell}}
\nc{\id}{\text{id}}
\nc\X{\mathcal X}

\nc{\Hom}{\on{Hom}}
\nc{\Lie}{\on{Lie}}
\nc{\Loc}{\on{Loc}}
\nc{\Pic}{\on{Pic}}
\nc{\Bun}{\on{Bun}}
\nc{\IC}{\on{IC}}
\nc{\Aut}{\on{Aut}}
\nc{\rk}{\on{rk}}
\nc{\Sh}{\on{Sh}}
\nc{\Perv}{\on{Perv}}
\nc{\pos}{{\on{pos}}}
\nc{\Conv}{\on{Conv}}
\nc{\Sph}{\on{Sph}}
\nc{\Sym}{\on{Sym}}
\nc{\BunBb}{\overline{\Bun}_B}
\nc{\Buno}{\overset{o}{\Bun}}
\nc{\BunPb}{{\overline{\Bun}_P}}
\nc{\BunBM}{\overline{\Bun}_{B(M)}}
\nc{\BunPbw}{{\widetilde{\Bun}_P}}
\nc{\BunBP}{\widetilde{\Bun}_{B,P}}
\nc{\GUb}{\overline{G/U}}
\nc{\GUPb}{\overline{G/U(P)}}

\nc{\Hhom}{\underline{\on{Hom}}}
\nc\syminfty{\on{Sym}^{\infty}}
\nc\lal{\ol{\lambda}}
\nc\xl{\ol{x}}
\nc\thl{\ol{\theta}}
\nc\nul{\ol{\nu}}
\nc\mul{\ol{\mu}}
\nc\Sum\Sigma
\nc{\oX}{\overset{o}{X}{}}

\nc{\M}{{\mathcal M}}
\nc{\N}{{\mathcal N}}
\nc{\F}{{\mathcal F}}
\nc{\D}{{\mathcal D}}
\nc{\Q}{{\mathcal Q}}
\nc{\Y}{{\mathcal Y}}
\nc{\G}{{\mathcal G}}
\nc{\E}{{\mathcal E}}
\nc{\CalC}{{\mathcal C}}
\nc\Dh{\widehat{\D}}
\renewcommand{\O}{{\mathcal O}}
\nc{\C}{{\mathcal C}}
\nc{\K}{{\mathcal K}}
\renewcommand{\H}{{\mathcal H}}
\renewcommand{\S}{{\mathcal S}}
\nc{\T}{{\mathcal T}}
\nc{\V}{{\mathcal V}}
\renc{\P}{{\mathcal P}}
\nc{\A}{{\mathcal A}}
\nc{\B}{{\mathcal B}}
\nc{\U}{{\mathcal U}}
\renewcommand{\L}{{\mathcal L}}
\nc{\Gr}{\on{Gr}}

\nc{\fA}{{\mathfrak A}}
\nc{\fP}{{\mathfrak P}}
\nc{\frn}{{\check{\mathfrak u}(P)}}
\nc{\p}{\overline{\mathfrak p}}
\nc{\q}{\overline{\mathfrak q}}
\nc\f{{\mathfrak f}}

\nc{\qo}{{\mathfrak q}}
\nc{\po}{{\mathfrak p}}
\nc{\s}{{\mathfrak s}}
\nc\w{\text{w}}
\renewcommand{\r}{{\mathfrak r}}

\newcommand{\tf}{{\mathfrak t}}

\nc\mathi\iota
\nc\Spec{\on{Spec}}
\nc\Mod{\on{Mod}}
\nc{\tw}{\widetilde{\mathfrak t}}
\nc{\pw}{\widetilde{\mathfrak p}}
\nc{\qw}{\widetilde{\mathfrak q}}
\nc{\jw}{\widetilde j}

\nc{\grb}{\overline{\Gr}}
\nc{\I}{\mathcal I}
\renewcommand{\i}{\mathfrak i}
\renewcommand{\j}{\mathfrak j}

\nc{\lambdach}{\check\lambda}
\nc{\Lambdach}{\check\Lambda}
\nc{\much}{\check\mu}
\nc{\omegach}{\check\omega}
\nc{\nuch}{\check\nu}
\nc{\etach}{\check\eta}
\nc{\alphach}{\check\alpha}
\nc{\betach}{\check\beta}
\nc{\rhoch}{\check\rho}

\nc{\Hb}{\overline{\H}}

\begin{document}

\title{Intersection cohomology of Drinfeld's compactifications}

%\author{K.$\Pi$.~ KABA$\Phi$H\,$\Sigma$}

\author{A.~Braverman, M.~Finkelberg, D.~Gaitsgory and I.~Mirkovi\'c}

\address{{\it Address}:\newline
A.B. and D.G.: Dept. of Math., Harvard University,
Cambridge MA 02138, USA; \newline
M.F.: Indep. Moscow Univ., 11 Bolshoj Vlasjevskij per., 
Moscow 121002, Russia; \newline
I.M.: Dept. of Math. and Stat.,
U-Mass. at Amherst,
Amherst MA 01003, USA.}

\email{\newline
braval@math.harvard.edu; fnklberg@mccme.ru;
gaitsgde@math.harvard.edu; \newline
mirkovic@math.umass.edu}

\maketitle

%\begin{flash}
%{\small \But, a ruined man, it's not my fault \\
%I'm only trying, poor devil, to make ends meet \\
%The almighty gods ought to have taken the trouble \\
%to create a fourth, an honest man. \\
%I would gladly have gone along with him. \\
%K.$\Pi$.~ KABA$\Phi$H\,$\Sigma$}
%\end{flash}

\section*{Introduction} \label{intr}

\ssec{}

This paper resulted by merging several points of view on a remarkable 
geometric object introduced by V.~Drinfeld. In the incarnation studied 
here, this is a relative compactification of the moduli stack of principal 
bundles with respect to a parabolic subgroup (of a given reductive group)
over a curve.

Since its discovery several years ago, Drinfeld's compactification has
found several remarkable applications in the geometric representation
theory, some of which are discussed in this Introduction. These applications
include a construction of the ``correct'' geometric Eisenstein series 
functor (cf. \cite{BG}); a geometric study of quantum groups at a root 
of unity and representations of Lie algebras in positive characteristic;
a realization of the combinatorial pattern introduced by Lusztig in \cite{Lu1} 
in terms of intersection cohomology. In particular, one obtains 
an interesting and unexpected relation between Eisenstein series and semi-infinite 
cohomology of quantum groups at a root of unity (cf. \cite{FFKM}).

The main result of this paper is the description of the intersection cohomology 
sheaf of Drinfeld's compactification. 

\ssec{The space $\BunPbw$}
Let $X$ be a smooth complete curve, $G$ a reductive group and 
$P\subset G$ a parabolic subgroup.
Let us denote by $\Bun_G$ the moduli stack of principal $G$-bundles on 
$X$ and by $\Bun_P$ the moduli stack of $P$-bundles. 

The inclusion of $P$ in $G$ gives rise to a representable morphism 
${\mathfrak p}:\Bun_P\to \Bun_G$, and it is a problem that arises most
naturally in geometric representation theory to look for 
a relative compactification of $\Bun_P$ along the fibers of this map. 

A construction of such a compactification was suggested by Drinfeld.
In this way we obtain a new algebraic stack, denoted 
$\BunPbw$ endowed with a map $\pw:\BunPbw\to\Bun_G$ 
such that $\pw$ is now {\it proper}, and $\Bun_P$ is contained
inside $\BunPbw$  as an open substack.

The main complication, as well as the source of interest,
is the fact that the stack $\BunPbw$ is singular. 
Basically, the present paper is devoted to study of singularities
of $\BunPbw$.

\ssec{Eisenstein series}

The authors of \cite{BG} considered the following problem:
For $X$, $G$ and $P$ as before, let $M$ be the Levi factor of $P$ and $\Bun_M$
the corresponding moduli stack. Using $\Bun_P$ one can introduce a naive 
{\it Eisenstein series} functor, which maps the derived category
$\on{D}^b(\Bun_M)$ to $\on{D}^b(\Bun_G)$, and which is an analogue of 
the usual Eisenstein series operator in the theory of automorphic forms.

However, since the map $\po$ is not proper, this functor is 
``not quite the right one'' from the geometric point if view. 
For example, it does not commute with the Verdier duality and 
does not preserve purity, etc. It turns out that using the 
compactified stack $\BunPbw$, one can introduce another functor $\on{Eis}^G_M$ 
between the same categories, which will have much better properties. The authors 
of {\it loc.cit.} called it the geometric Eisenstein series functor.

One of the results announced (but not proved) in \cite{BG} was Theorem 2.2.12, 
which compared the naive Eisenstein series functor with the corrected one. 
Very roughly, the new functor $\on{Eis}^G_M$ is the product of the old functor
and a (geometric version of) the appropriate $L$-function.

The proof of this theorem essentially amounts to an explicit description of 
the {\it intersection cohomology} sheaf on $\BunPbw$, in terms of
the combinatorics of the Langlands dual Lie algebra $\check{\mathfrak g}$, 
or more precisely in terms of the symmetric algebra $\on{Sym}^\bullet(\frn)$ 
(here $\frn$ is the Lie algebra of the unipotent 
radical of the corresponding parabolic in $\check G$), viewed as an 
$\check M$-module.

In the present paper we establish the required explicit description of
$\IC_{\BunPbw}$. 

\ssec{Semi-infinite flag variety and quantum groups}

Let us now explain another perspective on the behavior of the above $\on{IC}$ sheaf.
Consider the semi-infinite flag ``variety'' $G((t))/B((t))$, where $B$ is the Borel subgroup
of $G$. Since the pioneering work of Feigin and Frenkel \cite{FF}, people were 
trying to develop the theory of perverse sheaves (constructible with respect to a given
stratification) on $G((t))/B((t))$ and, in particular, to compute the $\on{IC}$ sheaf on it. 

The problem is that $G((t))/B((t))$ is very essentially infinite-dimensional, so 
that the conventional theory of perverse sheaves, defined for schemes of finite type,
was not applicable. 
Since it was (and still is) not clear whether there exists a direct 
definition of perverse sheaves on $G((t))/B((t))$, the authors of 
\cite{FFKM} proposed the following solution.

\medskip

They introduced certain finite-dimensional varieties, called the 
Zastava spaces $Z^\mu$ in terms of maps of a projective line into the flag variety 
(the parameter $\mu$ is the degree of the map, i.e. it is an element in the 
coroot lattice of $G$). They argued that
the these spaces provided ``models'' for $G((t))/B((t))$ from the point of view
of singularities of various strata.

Moreover, it was shown in \cite{FFKM} that the stalks of certain perverse sheaves 
on $Z^\mu$ were given by the periodic Kazhdan-Lusztig polynomials 
of \cite{Lu1}, 
and this agrees with the anticipated answer for $G((t))/B((t))$. Therefore, the Zastava spaces
provide a geometric interpretation for Lusztig's combinatorics. This, combined with the
earlier work of Feigin and Frenkel allowed the authors of \cite{FFKM} to come up with a 
conjecture that ties a certain category of perverse sheaves on $G((t))/B((t))$, thought of 
as sheaves on $Z^\mu$, with the category of representations of the small quantum group.

\medskip

The key fact for us is that the Zastava spaces $Z^\mu$, appropriately 
generalized in the case of an arbitrary parabolic $P$, provide ``local models'' 
for the singularities of the stack $\BunPbw$ as well. More precisely, we show
that the parabolic Zastava spaces and the stack $\BunPbw$ are, locally in the 
smooth topology, isomorphic to one another.

Here we are dealing with the following remarkable phenomenon: the stack $\BunPbw$
is defined via the global curve $X$ and it classifies $P$-bundles on our curve which
have degenerations at finitely many points. Therefore, one may wonder whether
the singularities of $\BunPbw$ near a point corresponding to such a degenerate 
$P$-bundle would depend only on what is happening at the points of degeneration, 
and this is what our comparison with the Zastava spaces actually proves.

\ssec{Intersection cohomology on parabolic Zastava spaces}

Thus, our problem is reduced to the computation of the $\on{IC}$ sheaf on the parabolic 
Zastava spaces $Z^\theta$ ($\theta$ now is an element of a certain quotient lattice).

To carry out the computation we employ an inductive procedure on the parameter $\theta$.
(We should note that the present argument is quite different and is in fact simpler than
the one used in \cite{FFKM} to treat the case of $P=B$.)

The basic characteristic of the Zastava spaces, discovered in \cite{FFKM}, is that
they are local in nature, which expresses itself in the {\it Factorization property}. 
Namely, Each $Z^\theta$ is fibered over the configuration space 
$X^\theta$, equal to the product of the corresponding
symmetric powers of $X$, i.e. $X^\theta=\underset{\i}\Pi\, X^{(n_\i)}$. 
If $\theta=\theta_1+\theta_2$
there is an isomorphism
$ Z^\theta\simeq 
\ Z^{\theta_1}\times Z^{\theta_2}$,\
{\em after both sides are restricted}
to the subspace
$(X^{\theta_1}\times X^{\theta_2})_{disj}$ of 
$X^{\theta_1}\times X^{\theta_2}$ corresponding
to non-intersecting configurations.

Moreover, the fiber of $Z^\theta$ over a given point of $X^\theta$
is a product of intersections of {\it semi-infinite orbits} in the 
{\it affine Grassmannian} $\Gr_G$ corresponding
to $G$. The whole space $Z^\theta$ can be thought of as a 
subspace of the corresponding version of the 
Beilinson-Drinfeld affine Grassmannian. 

\smallskip

This interpretation via the affine Grassmannian explains the connection between the
stalks of $\on{IC}_{Z^\theta}$ and the dual Lie algebra $\check{\mathfrak g}$:
the link is provided by the Drinfeld-Ginzburg-Lusztig-Mirkovi\'c-Vilonen theory 
of spherical perverse sheaves on $\Gr_G$ and their connection to $\on{Rep}(\check G)$, 
cf. \cite{MV}.

\ssec{The naive compactification}

To conclude, let us mention that in addition to $\BunPbw$, the stack $\Bun_P$ of 
parabolic bundles admits another, in a certain sense more naive, relative 
compactification, which we denote $\BunPb$. This second compactification
can be though of as a blow-down of $\BunPbw$; in particular, we have a proper
representable map $\r:\BunPbw\to \BunPb$.

By the Decomposition theorem, $\r_!(\IC_{\BunPbw})$ contains $\IC_{\BunPb}$ 
as a direct summand. In the last section
we give an explicit description of both $\r_!(\IC_{\BunPbw})$ and $\IC_{\BunPb}$. 
The answer for the stalks of $\on{IC}_\BunPb$ is formulated in terms of 
$\on{Sym}^\bullet(\frn^f)$ where $f\in\Lie\check M$ is a principal nilpotent.
 
Note that the latter vector space is exactly the one appearing in ~\cite{W}.
Therefore, the stack $\BunPb$ provides a geometric object whose singularities
reproduce the parabolic version of the periodic Kazhdan-Lusztig polynomials of ~\cite{W}.

\ssec{Notation}

In this paper we will work over the ground field $\Fq$. However, the reader will readily check
that all our results can be automatically carried over to the characteristic $0$ situation.

\medskip

Throughout the paper, $G$ will be a connected reductive group over $\Fq$, assumed
split. Moreover, we will assume that its derived group $G'=[G,G]$ is simply connected.

We will fix a Borel subgroup $B\subset G$ and let $T$ be the 
``abstract'' Cartan,
i.e., $T=B/U$, where $U$ is the unipotent radical of $B$. We will 
denote by $\Lambda$ the {\it coweight} lattice of
$T$ and by $\check \Lambda$--its dual, i.e. the weight lattice; $\langle\cdot,\cdot\rangle$ is the canonical
pairing between the two.

The semi--group of dominant coweights
(resp., weights) will be denoted by $\Lambda_G^+$ (resp., by $\Lambdach_G^+$). The set of vertices
of the Dynkin diagram of $G$ will be denoted by $\I$; for each $\i\in\I$ there corresponds a simple
coroot $\alpha_i$ and a simple root $\alphach_i$. The set of positive coroots will be denoted by $\Delta$
and their positive span inside $\Lambda$, by $\Lambda_G^{\on{pos}}$. (Note that, since $G$ is not semi--simple,
$\Lambda_G^+$ is not necessarily contained in $\Lambda_G^{\on{pos}}$.) By $\rhoch\in\Lambdach$ we will denote 
the half sum of positive roots of $G$ and by $w_0$ the longest element in the Weyl group of $G$.
For $\lambda_1,\lambda_2\in\Lambda$ we will write that 
$\lambda_1\geq \lambda_2$ if $\lambda_1-\lambda_2\in\Lambda_G^{\on{pos}}$,
and similarly for $\Lambdach_G^+$.

Let $P$ be a standard proper parabolic of $G$, i.e. $P\supset B$; let $U(P)$ be its unipotent radical
and $M:=P/U(P)$--the Levi factor. To $M$
there corresponds a sub-diagram $\I_M$ in $\I$. We will denote by $\Lambda_M^+\supset \Lambdach_G^+$, 
$\Lambda_M^{\on{pos}}\subset \Lambda_G^{\on{pos}}$, $\rhoch_M\in\Lambdach$, $w_0^M\in W$, $\underset{M}\geq$ etc.
the corresponding objects for $M$.

To a dominant coweight $\lambdach\in\Lambdach$ one attaches the Weyl 
$G$-module, denoted $\V^{\lambdach}$, with a fixed highest weight
vector $v^{\lambdach}\in \V^{\lambdach}$. For a pair of weights 
$\lambdach_1,\lambdach_2$, there is a canonical map
$\V^{\lambdach_1+\lambdach_2}\to \V^{\lambdach_1}\otimes \V^{\lambdach_2}$ that sends 
$v^{\lambdach_1+\lambdach_2}$ to $v^{\lambdach_1}\otimes v^{\lambdach_2}$.

Similarly, for $\nuch\in \Lambdach_M^+$, we will denote by $\U^{\nuch}$ the corresponding $M$--module. 
There is a natural functor $\V\mapsto \V^{U(P)}$ from the category of $G$--modules to that of $M$--modules and 
we have a canonical morphism $\U^{\lambdach}\to (\V^{\lambdach})^{U(P)}$, which sends the corresponding
highest weight vectors to one another.

\medskip

Unless specified otherwise, we will work with the perverse t-structure on the category
of constructible complexes over various schemes and stacks. If $\S$ is a constructible complex,
$h^i(\S)$ will denote its $i$-th perverse cohomology sheaf. The intersection
cohomology sheaves are normalized so that they are pure of weight $0$.
In other words, for a smooth variety $Y$ of dimension $n$, 
$\IC_Y\simeq (\Ql(\frac{1}{2})[1])^{\otimes n}$, where $\Ql(\frac{1}{2})$
corresponds to a chosen once and for all square root of $q$ in $\Ql$.

\ssec{Acknowledgments}

The first and fourth author acknowledge  partial support from the NSF.
A part of the work was completed while all authors 
were visiting  the Institute for Advanced Study in Princeton.
In the course of writing this paper, the second author enjoyed the 
hospitality and support of the Institut de Math\'ematiques de Luminy. 
In addition, he is grateful 
to A.~Kuznetsov for discussions of $\BunPb$.
The third author was supported by the Clay Mathematics Institute and 
by a grant from the NSF.

\section{Drinfeld's compactifications}

\ssec{}

Let $P$ be a parabolic subgroup of $G$ and let $H\subset P$ be either the unipotent radical $U(P)$ or $[P,P]$. 
Consider the following functor Schemes $\to$ Categories:

To a scheme $S$ we associate the category of triples $(\F_G,\F_{P/H},\kappa)$, where $\F_G$ (resp., $\F_{P/H}$) is a principal $G$-bundle
(resp., a principal $P/H$-bundle) on $X\times S$ and $\kappa$ is a collection of maps of coherent sheaves 
$$\kappa^\V:(\V^H)_{\F_{P/H}}\to \V_{\F_G},$$
defined for every $G$-module $\V$,
such that for every geometric point $s\in S$, $\kappa^\V_s$ is injective, and such that
the Pl\"ucker relations hold. This means that for $\V$ being the trivial representation, $\kappa^\V$
is the identity map $\O_{X\times S}\to \O_{X\times S}$ and for a morphism $\V_1\otimes \V_2\to \V_3$, the diagram
$$
\CD
(\V_1^H)_{\F_{P/H}}\otimes (\V_2^H)_{\F_{P/H}} @>{\kappa^{\V_1}\otimes \kappa^{\V_2}}>>
(\V_1\otimes \V_2)_{\F_G}  \\
@VVV    @VVV  \\
(\V_3^H)_{\F_{P/H}} @>{\kappa^{\V_3}}>> (\V_3)_{\F_G}
\endCD
$$
commutes.

Note that the data of $\kappa$ can be reformulated differently, 
using Frobenius reciprocity: for a $P/H$-module $\U$,
let $\on{Ind}(\U)$ denote the corresponding induced $G$-module, 
i.e. $\Hom_{P/H}(\U,\V^H)\simeq \Hom_G(\on{Ind}(\U),\V)$
for a $G$-module $\V$, functorially in $\V$. 
Then the data of $\kappa$ is the same as a collection of maps $\kappa$ defined
for every $P/H$-module $\U$:
$$\kappa^{\U}:\U_{\F_{P/H}}\to (\on{Ind}(\U))_{\F_G},$$
which satisfy the Pl\"ucker relations 
in the sense that this map is again the identity for the trivial representation and
for $\U_3\to \U_1\otimes \U_2$ the diagram
$$
\CD
(\U_3)_{\F_{P/H}}   @>>> (\on{Ind}(\U_3))_{\F_G}  \\
@VVV      @VVV      \\
(\U_1\otimes \U_2)_{\F_{P/H}}  @>>>  (\on{Ind}(\U_1\otimes \U_2))_{\F_G}
\endCD
$$
commutes.

In particular, for $H=[P,P]$, it is sufficient to specify the value of $\kappa$ on
$1$-dimensional representations of $P/H$, since this group is commutative.

\medskip

For a fixed $S$, it is clear that triples $(\F_G,\F_{P/H},\kappa)$ form a groupoid, and for a map $S_1\to S_2$
there is a natural (pull-back) functor between the corresponding groupoids. In addition, there is a natural forgetful morphism
from this functor to the functor represented by the stack $\Bun_G$: out of $(\F_G,\F_{P/H},\kappa)$ we "remember"
only $\F_G$, which is a $G$-bundle on $X\times S$.

The following facts are proven in \cite{BG}:

\begin{thm}
For both choices of $H$ the above functor $S\mapsto (\F_G,\F_{P/H},\kappa)$ is representable 
by an algebraic stack, such that its map to $\Bun_G$ is representable and proper.
\end{thm} 

We denote the corresponding stacks by $\BunPbw$ (when $H=U(P)$) and by $\BunPb$ (when $H=[P,P]$).
Their projections to $\Bun_G$ will be denoted by $\pw$ and $\p$, respectively.

\medskip

Note that we have a natural map from the stack $\Bun_P$ to both $\BunPbw$ and $\BunPb$. Indeed, 
in both cases, a $P$-bundle on $X\times S$ is the same as a triple 
$(\F_G,\F_{P/H},\kappa)$, for which the maps $\kappa^\V$ are injective bundle maps. 

Since the condition for a map between vector bundles to be {\it maximal} (maximal means to be an injective bundle map)
is open, the above maps $\Bun_P\to \BunPbw$ and $\Bun_P\to \BunPbw$ are open embeddings. The following 
is also established in {\it loc. cit}:

\begin{thm}
$\Bun_P$ is dense in both $\BunPbw$ and $\BunPb$.
\end{thm}

Finally, note that since $U(P)$ is contained in $[P,P]$ we have a natural map
$\r:\BunPbw\to\BunPb$, which is proper and whose restriction to $\Bun_P$ is the identity map.

\ssec{}

By construction, a point of $\BunPbw$ (with values in a field) defines a $P$-bundle over an open subset $X^0$
of the curve $X$ (in fact $X^0$ is precisely the locus, where the maps $\kappa^\V$ are maximal embeddings).
We will now describe the partition of these stacks according to the behavior of $(\F_G,\F_{P/H},\kappa)$
on $X-X^0$. First, we will treat the case of $\BunPb$.

Let $M$ be the Levi factor of $P$. We choose a splitting $M\hookrightarrow P$; in particular we
have a well-defined opposite parabolic $P^-$ such that $P\cap P^-=M$. We will denote by $\I_M$ the corresponding
Dynkin sub-diagram of $\I$.

Let us denote by $\Lambda_{G,P}$ the lattice of cocharacters of the torus $P/H\simeq M/[M,M]$. We have
the natural projection $\Lambda\to \Lambda_{G,P}$, whose kernel is the span of $\alpha_\i$, $\i\in\I_M$. We
will denote by $\Lambda_{G,P}^{\on{pos}}$ the sub-semigroup of $\Lambda_{G,P}$ spanned by the images of $\alpha_\j$,
$\j\in\I-\I_M$.

Given $\theta\in \Lambda_{G,P}^{\on{pos}}$, we will denote by ${\mathfrak A}(\theta)$ the elements of the set
of decompositions of $\theta$ as a sum of non-zero elements of $\Lambda_{G,P}^{\on{pos}}$. In other words,
each ${\mathfrak A}(\theta)$ is a way to write $\theta=\underset{k}\Sum\, n_k\cdot \theta_k$, where
$\theta_k$'s belong to $\Lambda_{G,P}^{\on{pos}}-0$ and are pairwise distinct.

For ${\mathfrak A}(\theta)$ we will denote by $X^{{\mathfrak A}(\theta)}$ the corresponding partially symmetrized
power of the curve, i.e. $X^{{\mathfrak A}(\theta)}=\underset{k}\Pi\, X^{(n_k)}$. We will denote by 
$\overset{o}X{}^{{\mathfrak A}(\theta)}$ the complement of the diagonal divisor in $X^{{\mathfrak A}(\theta)}$.

\begin{prop}  \label{naivestrata}
There exists a locally closed embedding $j_{{\mathfrak A}(\theta)}:\overset{o}X{}^{{\mathfrak A}(\theta)}\times 
\Bun_P\to \BunPb$. Every field-valued point of $\BunPb$ belongs to the image a unique $j_{{\mathfrak A}(\theta)}$.
\end{prop}

Although the proof is given in \cite{BG}, let us indicate the construction of $j_{{\mathfrak A}(\theta)}$.

\begin{proof}

Let $\F_P$ be a point of $\Bun_P$ and $x^{{\mathfrak A}(\theta)}$ be a point of $\overset{o}X{}^{{\mathfrak A}(\theta)}$, i.e.
$x^{{\mathfrak A}(\theta)}$ consists of a collection of pairwise distinct points $x_1,...,x_n$ of $X$, to each of 
which there is an assigned element
$\theta_k\in \Lambda_{G,P}^{\on{pos}}$. We define the corresponding point of $\BunPb$ as follows:

In the triple $(\F_G,\F_{M/[M,M]},\kappa)$, $\F_G$ is the $G$-bundle induced from $\F_P$. Let now
$\F'_{M/[M,M]}$ be the $M/[M,M]$-bundle induced from $\F_P$. For each $G$-dominant weight $\lambdach$ orthogonal
to $\on{Span}(\alpha_\i)$, $\i\in\I_M$, let us denote by $\L^{\lambdach}_{\F'_{M/[M,M]}}$ the corresponding
associated line bundle. By construction, we have the injective bundle maps
$$\kappa'{}^{\lambdach}:\L^{\lambdach}_{\F'_{M/[M,M]}} \to \V^{\lambdach}_{\F_G}$$
(here $\V^{\lambdach}$ is the Weyl module corresponding to $\lambdach$), which satisfy Pl\"ucker relations.

We set $\F_{M/[M,M]}:=\F'_{M/[M,M]}(-\Sum \, \theta_k\cdot x_k)$. The corresponding line bundles
$\L^{\lambdach}_{\F_{M/[M,M]}}$ are then $\L^{\lambdach}_{\F'_{M/[M,M]}}(-\Sum\, \langle \theta_k,\lambdach\rangle\cdot x_k)$.
Thus, by composing we obtain the maps
$$\kappa^{\lambdach}:\L^{\lambdach}_{\F_{M/[M,M]}} \to \V^{\lambdach}_{\F_G},$$
which are easily seen to satisfy the Pl\"ucker relations, as required.

\end{proof}

Let us denote the image of $j_{{\mathfrak A}(\theta)}$ in $\BunPb$ by $_{{\mathfrak A}(\theta)}\BunPb$. 
It is easy to see that the union of $_{{\mathfrak A}(\theta)}\BunPb$ is also a locally closed sub-stack of $\BunPb$,
which we will denote by $_\theta\BunPb$.

If $\theta$ is the projection of $\underset{i\in \I-\I_M}\Sigma\,  n_\i\cdot \alpha_\i$, set 
$X^\theta=\underset{\i\in\I-\I_M}\Pi\, X^{(n_\i)}$. By definition, $X^\theta$ is stratified by the spaces
$\overset{o}X{}^{{\mathfrak A}(\theta)}$ for all possible ${\mathfrak A}(\theta)$. As in \propref{naivestrata},
we have a locally closed embedding $X^\theta\times \Bun_P\to\BunPb$, whose image is our $_\theta\BunPb$.

\medskip

Let us denote by $_{{\mathfrak A}(\theta)}\BunPbw\subset \BunPbw$ the preimage of $_{{\mathfrak A}(\theta)}\BunPb$
under the map $\r$. Our next goal is to give a more explicit description of each $_{{\mathfrak A}(\theta)}\BunPbw$.

\ssec{}

Let $x\in X$ be a  point. Recall that the affine Grassmannian $\Gr_G$ is the ind-scheme representing the functor
Schemes $\to$ Sets that attaches to a scheme $S$ the set of pairs $(\F_G,\beta)$, where $\F_G$ is a principal $G$-bundle
on $X\times S$ and $\beta$ is an isomorphism $\F_G|_{(X-x)\times S}\simeq \F^0_G|_{(X-x)\times S}$, where $\F^0_G$
is {\it the trivial} $G$-bundle. Sometimes, in order to emphasize the dependence on the point $x$, we will write $\Gr_{G,x}$. 
Note that by letting 
$x$ move along the curve, we obtain a relative version of $\Gr_G$, which will be denoted $\Gr_{G,X}$.

In the same way one defines the affine Grassmannians for the groups $M$, $P$, $P^-$.

For every $G$-dominant coweight $\lambda$ one defines a (finite-dimensional) closed
sub-scheme $\grb^\lambda_G\subset \Gr_G$ by the condition that $(\F_G,\beta)\in \grb^\lambda_G$ if for every $G$-module
whose weights are $\leq \lambdach$, the meromorphic map 
$\beta^\V:\V_{\F_G}\to \V_{\F^0_G}\simeq \V\otimes \O_{X\times S}$ has
a pole of order $\leq \langle \lambda,-w_0(\lambdach)\rangle$ along $x\times S$.

\medskip

Now let $\theta$ be an element of $\Lambda_{G,P}^{\on{pos}}$. We define the element $\flat(\theta)\in \Lambda$
as follows: if $\theta$ is the projection under $\Lambda\to \Lambda_{G,P}$ of $\widetilde{\theta}\in \on{Span}(\alpha_\j)$, 
$\j\in \I-\I_M$, then $\flat(\theta)=w_0^M(\widetilde{\theta})$, where $w_0^M$ is the longest element in the Weyl
group of $M$. Note that by construction, 
$\flat(\theta)$ is $M$-dominant; in particular, it makes sense to consider
$\grb_M^{\flat(\theta)}$. 

Consider the functor that attaches to a scheme $S$ the set of pairs 
$(\F_M,\beta^\V)$, where $\beta^\V$ is 
an embedding of coherent sheaves defined for every $G$-module $\V$:
$$\beta_M^\V:\V^{U(P)}_{\F_M}\hookrightarrow \V^{U(P)}_{\F^0_M}$$ such that

\smallskip

\noindent 1) The Pl\"ucker relations hold in the same sense as before.

\smallskip

\noindent 2) If $\V^{U(P)}$ is 1-dimensional corresponding to the character $\nuch$ of $M$, then $\beta_M^\V$ identifies
$\V^{U(P)}_{\F_M}:=\L^{\nuch}_{\F_M}$ with $\L^{\nuch}_{\F^0_M}(-\langle \theta,\nuch\rangle \cdot x) \simeq\O_{X\times S}
(-\langle \theta,\nuch\rangle \cdot x)$. (In the last formula we have used the fact that $\nuch$ and $\theta$ belong to
mutually dual lattices.)

\medskip

The following proposition is proved in \cite{BG}, but we will sketch the argument due to its importance:

\begin{prop}  \label{pospart}
The above functor is representable by a finite dimensional closed sub-scheme, denoted $\Gr_M^{+,\theta}$, 
of $\Gr_M$. We have an inclusion $\grb_M^{\flat(\theta)}\hookrightarrow \Gr_M^{+,\theta}$, which
induces an isomorphism on the level of reduced schemes.
\end{prop}

\begin{proof}

Let $(\F_M,\beta_M^\V)$ be an $S$-point of $\Gr_M^{+,\theta}$. To construct a map of functors,
$\Gr_M^{+,\theta}\to\Gr_M$ we must exhibit the maps $\beta^\U_M:\U_{\F_M}\to \U_{\F^0_M}$
for all $M$-modules $\U$ and not just for those of the form $\V^{U(P)}$. However, we can do that because
any $\U$ can be tensored with a $1$-dimensional representation of $M$ corresponding to a $G$-dominant weight $\nuch$,
so that the new representation will be of the form $\V^{U(P)}$.

By construction, the above map $\Gr_M^{+,\theta}\to\Gr_M$ is a closed embedding. The fact that
$\Gr_M^{+,\theta}$ is a scheme (and not an ind-scheme) follows from the fact that we can choose $\V$ such that
$\V^{U(P)}$ is faithful as a representation of $M$.

Let $(\F_M,\beta_M)$ be an $S$-point of $\grb_M^{\flat(\theta)}$. Then if $\V$ is a $G$-module (whose weights, we can
suppose, are $\leq \lambdach$ for some $G$-dominant weight $\lambdach$), the maps 
$\beta_M^{\V^{U(P)}}:\V^{U(P)}_{\F_G}\to \V^{U(P)}_{\F^0_G}$ are regular, since 
$\langle \flat(\theta),-w^M_0(\lambdach)\rangle\leq 0$, by the definition of $\flat(\theta)$.

Hence, $\grb_M^{\flat(\theta)}$ is contained in $\Gr_M^{+,\theta}$. To show that this inclusion is
an isomorphism on the level of reduced schemes, one has to check that the fact that $\grb_M^\lambda\subset \Gr_M^{+,\theta}$
implies that $\flat(\theta)-\lambda$ is a sum of positive coroots of $M$, which is obvious.

\end{proof}

\ssec{}

Now let us consider the following relative version of the above situation. Let $\theta$ be as above and given an element
${\mathfrak A}(\theta)$ let us consider the space $\Gr_M^{+,{\mathfrak A}(\theta)}$ over 
$\overset{o}X{}^{{\mathfrak A}(\theta)}$,
whose fiber of $x^{{\mathfrak A}(\theta)}=\Sum\, \theta_k\cdot x_k\in \overset{o}X{}^{{\mathfrak A}(\theta)}$ equals
$\underset{k}\Pi\, \Gr_{M,x_k}^{+,\theta_k}$.

In addition, we can generalize this further, by replacing the trivial $M$-bundle in the definition of $\Gr_M$ 
by an arbitrary background $M$-bundle $\F_M^b$. By letting $\F^b_M$ vary along the universal family, i.e. $\Bun_M$,
we obtain the relative version $\H_{M}$ of $\Gr_{M,X}$, which is fibered over $X\times\Bun_M$ and is known in the literature as the
{\it Hecke stack}. The relative version of $\Gr_M^{+,{\mathfrak A}(\theta)}$ will be denoted by $\H_M^{+,{\mathfrak A}(\theta)}$ and
it is by definition fibered over $\overset{o}X{}^{{\mathfrak A}(\theta)}\times\Bun_M$. 

\begin{prop}   \label{descrBunPbw}
There exists a canonical isomorphism 
$$_{{\mathfrak A}(\theta)}\BunPbw\simeq \Bun_P\underset{\Bun_M}\times \H_M^{+,{\mathfrak A}(\theta)},$$
such that the projection $\r$ onto $_{{\mathfrak A}(\theta)}\BunPb$ on the LHS corresponds to the natural map of the RHS
to $\Bun_P\times \overset{o}X{}^{{\mathfrak A}(\theta)}$.
\end{prop}

The proof is given in \cite{BG} and is, in fact, an easy consequence of \propref{pospart} above.

\ssec{}

Finally, we are able to state \thmref{restr}, which is the main result of this paper. 

\medskip

First, let us recall the category of spherical perverse sheaves on $\Gr_G$, which by definition consists
of direct sums of perverse sheaves $\IC_{\grb_G^\lambda}$, as $\lambda$ ranges over the set of $G$-dominant coweights.
It is known that this category possesses a tensor structure, given by the convolution product, and as a tensor
category it is equivalent to the category $\on{Rep}(\check G)$ of finite-dimensional representations of the Langlands dual
group $\check G$. In particular, we have the functor $\on{Loc}:\on{Rep}(\check G)\to\on{Prev}(\Gr_G)$, such that
the irreducible representation of $\check G$ with h.w. $\lambda$ goes over under this functor to $\IC_{\grb_G^\lambda}$.

We will use the above definitions for $M$, rather than for $G$. Recall that $\Lambda_{G,P}$ can be canonically identified
with the lattice of characters of $Z(\check M)$; for $\theta\in\Lambda_{G,P}$ and an $\check M$-representation $V$, we will
denote by $V_\theta$ the direct summand of $V$ on which $Z(\check M)$ acts according to $\theta$.

Recall that the nilpotent radical of the dual parabolic 
$\check{\mathfrak u}(P)$ is naturally a representation of the group $\check M$ and let us observe
that for $\theta\in \Lambda_{G,P}^{\on{pos}}$, the sub-representation $\check{\mathfrak u}(P)_\theta$ is irreducible.

\begin{lem}
For $\theta\in \Lambda_{G,P}^{\on{pos}}$, the perverse sheaf $\on{Loc}(\on{Sym}(\check{\mathfrak u}(P))_\theta)$
is supported on $\grb_M^{\flat(\theta)}$.
\end{lem}

\medskip

Now, let us fix the notation for the relative versions of the functor $\on{Loc}$. First, we will denote by $\on{Loc}_X$
the corresponding functor $\on{Loc}:\on{Rep}(\check M)\to\on{Perv}(\Gr_{M,X})$. Secondly, given $\theta$ and an
element ${\mathfrak A}(\theta)$, corresponding 
to $\theta=\underset{k}\Sum\, n_k\cdot \theta_k$, and an $\check M$-representation
$V$, we will denote by $\on{Loc}_X^{{\mathfrak A}(\theta)}(V)$ the perverse sheaf on $\Gr_M^{+,{\mathfrak A}(\theta)}$,
whose fiber over $x^{{\mathfrak A}(\theta)}=\Sum\, \theta_k\cdot x_k\in \overset{o}X{}^{{\mathfrak A}(\theta)}$, i.e. over
$\underset{k}\Pi\, \Gr_{M,x_k}^{+,\theta_k}$, is
$$\underset{k}\boxtimes\, \on{Loc}(V_{\theta_k})\otimes (\Ql[1](\frac{1}{2}))^{\otimes |{\mathfrak A}(\theta)|},$$
where $|{\mathfrak A}(\theta)|=\underset{k}\Sum\, n_k$.

We will use the symbol $\on{Loc}_{\Bun_M,X}(V)$ (resp., $\on{Loc}_{\Bun_M,X}^{{\mathfrak A}(\theta)}(V)$)
to denote the corresponding perverse sheaves on $\H_{M}$ (resp., $\H_{M}^{+,{\mathfrak A}(\theta)}$).
Furthermore, by considering the tensor product $\Bun_P\underset{\Bun_M}\times \H_{M}$
we can define a perverse sheaf $\on{Loc}_{\Bun_P,X}(V)$ on it, which is $\IC_{\Bun_P}$ ``along the base''
(i.e. $\Bun_P$) and $\on{Loc}_{\Bun_M,X}(V)$ ``along the fiber'' (i.e. $\H_{M}$), and similarly for the perverse
sheaf $\on{Loc}_{\Bun_P,X}^{{\mathfrak A}(\theta)}(V)$ on $\Bun_P\underset{\Bun_M}\times \H_{M}^{+,{\mathfrak A}(\theta)}$.

\begin{thm} \label{restr}
The $*$-restriction of $\IC_{\BunPbw}$ to $_{{\mathfrak A}(\theta)}\BunPbw\simeq \Bun_P\underset{\Bun_M}\times 
\H_M^{{\mathfrak A}(\theta)}$
is isomorphic to
$$\on{Loc}_{\Bun_P,X}^{{\mathfrak A}(\theta)}(\underset{i\geq 0}\oplus\, \on{Sym}^i(\check{\mathfrak u}(P))\otimes \Ql(i)[2i])
\otimes (\Ql[1](\frac{1}{2}))^{\otimes -|{\mathfrak A}(\theta)|},$$
where $\underset{i\geq 0}\oplus\,\on{Sym}^i(\check{\mathfrak u}(P))\otimes \Ql(i)[2i]$ is viewed as a 
cohomologically graded $\check M$-module.
\end{thm}

\section{Zastava spaces}

Let $\theta$ be an element of $\Lambda_{G,P}^{\on{pos}}$. In this section 
we will introduce the Zastava spaces $Z^{\theta}$, which will be local models for $\BunPbw$.
 
\ssec{}

Let us recall the space $X^\theta$: if $\theta=\underset{i\in \I-\I_M}\Sigma\,  n_\i\cdot \alpha_\i$,
$X^\theta=\underset{\i\in\I-\I_M}\Pi\, X^{(n_\i)}$. One may alternatively view $X^\theta$ 
as the space classifying the following data:
$(\F_{M/[M,M]},\beta_{M/[M,M]})$, where $\F_{M/[M,M]}$ is a principal bundle with respect to the group $M/[M,M]$ on $X$
of degree $-\theta$, and $\beta_{M/[M,M]}$ is a system of embeddings defined for every $G$-dominant weight $\lambdach$
orthogonal to $\on{Span}(\alpha_\i)$, $\i\in\I_M$
$$\beta^{\lambdach}_{M/[M,M]}:\L^{\lambdach}_{\F_{M/[M,M]}}\to \L^{\lambdach}_{\F^0_{M/[M,M]}}\simeq \O_X,$$
such that $\beta^{\lambdach_1}\otimes \beta^{\lambdach_2}=\beta^{\lambdach_1+\lambdach_2}$.

Note that in the product $X\times X^\theta$ there is a natural incidence divisor, which we will denote by $\Gamma^\theta$.

Now, let us define the scheme $\Mod_M^{+,\theta}$.
By definition, its $S$-points are pairs $(\F_M,\beta_M)$, where $\F_M$ is an $M$-bundle on $X\times S$ 
such that the induced $M/[M,M]$-bundle is of degree $-\theta$, and $\beta_M$ is a system of embeddings of coherent sheaves
defined for every $G$-module $\V$
$$\beta^\V_M:(\V^{U(P)})_{\F_M}\hookrightarrow  (\V^{U(P)})_{\F^0_M},$$
such that for a pair of $G$-modules $\V_1$ and $\V_2$ we have a commutative diagram
$$
\CD
(\V_1^{U(P)})_{\F_M}\otimes (\V_2^{U(P)})_{\F_M}  @>{\beta^{\V_1}_M\otimes \beta^{\V_2}_M}>> 
(\V_1^{U(P)})_{\F^0_M}\otimes (\V_2^{U(P)})_{\F^0_M}  \\
@VVV   @VVV   \\
(\V_1\otimes \V_2)^{U(P)}_{\F_M}  @>{\beta_M^{\V_1\otimes\V_2}}>> (\V_1\otimes \V_2)^{U(P)}_{\F^0_M}.
\endCD
$$

It is easy to see, as in \propref{pospart}, that $\Mod_M^{+,\theta}$ is indeed representable by
a scheme of finite type.

By construction, we have a natural map $\pi_M:\Mod_M^{+,\theta}\to X^\theta$, which corresponds to taking
for $\V$ all possible $1$-dimensional $M$-modules. If $(\F_M,\beta_M)$ is an $S$-point of
$\Mod_M^{+,\theta}$, it follows as in \propref{pospart} that $\beta_M$ defines a trivialization of $\F_M$
on $X\times S-\Gamma_S^\theta$, where $\Gamma_S^\theta$ is the preimage of $\Gamma^\theta$ under
$X\times S\to X\times \Mod_M^{+,\theta}\to X\times X^\theta$. Moreover, we have:
$$\overset{o}X{}^{{\mathfrak A}(\theta)}\underset{X^\theta}\times \Mod_M^{+,\theta}\simeq \Gr_M^{+,{\mathfrak A}(\theta)}.$$
In particular, $\Mod_M^{+,\theta}|_{\Delta_X}\simeq \Gr_{M,X}^{+,\theta}$, where $\Delta_X\subset X^\theta$
is the main diagonal.

\ssec{}

Finally, we are ready to define $Z^\theta$. An $S$-point of $Z^\theta$ is a quadruple $$(\F_G,\F_M,\beta_M,\beta),$$
where $\F_G$ is a $G$-bundle on $X\times S$,
$(\F_M,\beta_M)$ is a point of $\Mod_M^{+,\theta}$ and $\beta$ is a trivialization of $\F_G$ on
$X\times S-\Gamma_S^\theta$, where $\Gamma^\theta_S$ is as above, such that the following two conditions are satisfied:

\begin{enumerate}

\item

For every $G$-module $\V$, the natural map $\V\to \V_{U(P^-)}$ extends to a regular surjective map of vector bundles
$$\V_{\F_G}\overset{\beta}\to \V_{\F^0_G}\to(\V_{U(P^-)})_{\F^0_M}\simeq \V_{U(P^-)}\otimes \O_{X\times S}.$$

\item

For every $G$-module $\V$, the natural map $\V^{U(P)}\to \V$ extends by means of $\beta$ and $\beta_M$
to a regular embedding of coherent sheaves
$$(\V^{U(P)})_{\F_M}\to \V_{\F_G}.$$

\end{enumerate}

\medskip

From the above definition, it follows that $Z^\theta$ is representable by an ind-scheme.
However, we will see later on that $Z^\theta$ is in fact a scheme, cf. \propref{identification}.

\medskip

We will denote by $\pi_P$ the natural map $Z^\theta\to \Mod_M^{+,\theta}$; by $\pi_G$ we will denote the composition
$\pi_M\circ \pi_P:Z^{\theta}\to X^\theta$.

\medskip

By definition, $Z^\theta$ contains as a sub-scheme the locus of those
$(\F_G,\F_M,\beta_M,\beta)$, for which the maps $(\V^{U(P)})_{\F_M}\to \V_{\F_G}$ are {\it maximal embeddings}, i.e. are
bundle maps. We will denote this sub-scheme by $Z_{max}^\theta$.

\smallskip

Observe now that there is a natural closed embedding: $\s^\theta:\Mod_M^{+,\theta}\to Z^\theta$. Indeed, to 
$(\F_M,\beta_M)\in \Mod_M^{+,\theta}$ we attach $(\F^0_G,\F_M,\beta_M,\beta^0)$, where $\beta^0$ is
the tautological trivialization of the trivial bundle. 

\medskip

\noindent{\it Remark.}
Note that for the definition of the Zastava space $Z^\theta$, the curve $X$ need not be complete.
Indeed, the only modification is the following:

In the definition of $\Mod_M^{+,\theta}$, instead of having pairs
$(\F_M,\beta_M)$ we can consider triples $(\F_M,\beta_M,S\to X^\theta)$
where $\beta_M$ is such that for every $1$-dimensional $M$-module $\V^{\lambdach}$
(such $\lambdach$ is automatically orthogonal to $\alpha_\i$, $\i\in\I_M$),
$\beta_M$ induces an isomorphism $\V^{\lambdach}_{\F_M}\simeq 
\O_{X\times S}(-\langle \theta,\lambdach\rangle\cdot \Gamma^\theta_S)$.

\ssec{Factorization property}

The fundamental property of the spaces $Z^\theta$ is their local behavior with respect to the base
$X^\theta$. 

Let $\theta=\theta_1+\theta_2$ with $\theta_i\in\Lambda_{G,P}^{\on{pos}}$ and let us denote by
$(X^{\theta_1}\times X^{\theta_2})_{disj}$ the open subset of the direct product $X^{\theta_1}\times X^{\theta_2}$,
which corresponds to $x^{\theta_1}\in X^{\theta_1}$, $x^{\theta_2}\in X^{\theta_2}$, such that the supports of
$x^{\theta_1}$ and $x^{\theta_2}$ are disjoint.

We have a natural \'etale map $(X^{\theta_1}\times X^{\theta_2})_{disj}\to X^\theta$.

\begin{prop}  \label{factor}
There is a  natural isomorphism:
$$(X^{\theta_1}\times X^{\theta_2})_{disj}\underset{X^\theta}\times Z^\theta\simeq 
(X^{\theta_1}\times X^{\theta_2})_{disj}\underset{X^{\theta_1}\times X^{\theta_2}} \times (Z^{\theta_1}\times Z^{\theta_2}).$$
\end{prop}

\begin{proof}

Let $x^{\theta_1}\times x^{\theta_2}$ be an $S$-point of $(X^{\theta_1}\times X^{\theta_2})_{disj}$. By 
definition, this means that the divisors $\Gamma^{\theta_1}_S$ and $\Gamma^{\theta_2}_S$ in $X\times S$
do not intersect. Let $(\F_G,\F_M,\beta_M,\beta)$ be an $S$-point of $Z^\theta$ which projects under $\pi_G$ to 
the corresponding point of $X^\theta$.

Set $(X\times S)^1=X\times S-\Gamma^{\theta_1}_S$, 
$(X\times S)^2=X\times S-\Gamma^{\theta_2}_S$, 
$(X\times S)^0=(X\times S)^1\cap (X\times S)^2$. By assumption, $(X\times S)^1\cup (X\times S)^2=X\times S$. We define a new
$G$-bundle $\F^1_G$ as follows: over $(X\times S)^1$, $\F^1_G$ is by definition the trivial 
bundle $\F^0_G$; over $(X\times S)^2$, $\F^1_G$ is identified with $\F_G$;
the data of $\beta$, being a trivialization of $\F_G$ over $(X\times S)^0$ defines a patching data for $\F^1_G$. By construction,
$\F^1_G$ is trivialized off $x^{\theta_1}$, let us denote this trivialization by $\beta^1$.

We introduce the second $G$-bundle $\F^2_G$ in a similar fashion: 
$\F_2|_{(X\times S)^2}\simeq \F_G|_{(X\times S)^2}$ and $\F_2|_{(X\times S)^1}\simeq \F^0_G|_{(X\times S)^1}$.
From the construction, $\F^2_G$ acquires a trivialization $\beta^2:\F^2_G|_{(X\times S)^2}\simeq \F^0_G|_{(X\times S)^2}$. 

In a similar way, from $(\F_M,\beta_M)$ we obtain two pairs $(\F^1_M,\beta^1_M)\in \Mod_M^{+,\theta_1}$ and 
$(\F^2_M,\beta^2_M)\in \Mod_M^{+,\theta_2}$, which project under $\pi_M$ to $x^{\theta_1}$ and $x^{\theta_2}$, respectively.

Thus, from the $S$-point $(\F_G,\F_M,\beta_M,\beta)$ we obtain two $S$-points $(\F^1_G,\F^1_M,\beta^1_M,\beta^1)$ and 
$(\F^2_G,\F^2_M,\beta^2_M,\beta^2)$ of
$Z^{\theta_1}$ and $Z^{\theta_2}$, respectively. The map in the opposite direction is constructed in the same way.

\end{proof}

In the course of the proof we have shown that the space $\Mod_M^{+,\theta}$ factorizes as well, i.e. we have
a natural isomorphism 
$$(X^{\theta_1}\times X^{\theta_2})_{disj}\underset{X^\theta}\times \Mod_M^{+,\theta}\simeq 
(X^{\theta_1}\times X^{\theta_2})_{disj}\underset{X^{\theta_1}\times X^{\theta_2}}\times 
\Mod_M^{+,\theta_1}\times \Mod_M^{+,\theta_2},$$ compatible with the factorization of $Z^\theta$. 
In addition, it is clear that the section $\s^\theta$ is compatible
with the factorizations in a natural way.

\ssec{The central fiber}

Consider the main diagonal $\Delta_X:X\to X^\theta$. For a fixed point
$x\in X$ let us consider the corresponding composition
$\Delta_x:\on{pt}\to X\to X^\theta$.

\medskip

The central fiber $\SS^\theta$ of $Z^\theta$ is by definition the preimage 
of the above point under $\pi_G:Z^\theta\to X^\theta$.
We will denote by $_0\SS^\theta$ the intersection $\SS^\theta\cap Z^\theta_{max}$.

\medskip

For $\theta\in\Lambda_{G,P}$, let $\Gr_P^\theta$ be the pre-image under 
$\Gr_P\to \Gr_{M/[M,M]}$ of the corresponding point-scheme
in $\Gr_{M/[M,M]}$. Both $\Gr_P^\theta$ and $\Gr_{U(P^-)}$ are locally 
closed sub-schemes of $\Gr_G$ and let us consider
their intersection $\Gr_P^\theta\cap \Gr_{U(P^-)}$. 

\begin{prop} \label{central fiber}
There is a natural identification $\Gr_P^\theta\cap \Gr_{U(P^-)}\simeq {}_0\SS^\theta$.
The map $_0\SS^\theta\overset{\pi_P}\to \Gr_M^{+,\theta}\hookrightarrow \Gr_M$
corresponds to $\Gr_P^\theta\cap \Gr_{U(P^-)}\hookrightarrow \Gr_P\to\Gr_M$.
\end{prop}

\begin{proof}

By construction, an $S$-point of $_0\SS^\theta$ is a data of a $G$-bundle $\F_G$ on $X\times S$, with
given reductions $\F_P$ and $\F_{P^-}$ to $P$ and $P^-$, respectively, such that
these reductions are mutually transversal on $(X-x)\times S$, the
$M$-bundle induced from $\F_{P^-}$ is trivialized, and the maps 
$$\L^{\lambdach}_{\F_{M/[M,M]}}\to \O_{X\times S}(-\langle \lambdach,\theta\rangle)$$
are isomorphisms for $G$-dominant characters $\lambdach$ of $M/[M,M]$.

Therefore, over $(X-x)\times S$ all the three principal bundles $\F_G$, $\F_P$ and $\F_{P^-}$
are trivialized in a compatible way, and the $M/[M,M]$-bundle induced from $\F_P$ is exactly
$\F^0_{M/[M,M]}(-\theta\cdot x)$. Thus, $(\F_P,\F_{P^-})$ indeed defines a point of 
$\Gr_P^\theta\cap \Gr_{U(P^-)}$. 

The map in the opposite direction, i.e. $\Gr_P^\theta\cap \Gr_{U(P^-)}\to {}_0\SS^\theta$,
is constructed in exactly the same way. The second assertion of the proposition, i.e.
the fact that $\pi_P$ becomes identified with $\Gr_P^\theta\cap \Gr_{U(P^-)}\hookrightarrow \Gr_P\to\Gr_M$
is obvious from the construction.

\end{proof}

As a consequence, we obtain that since $Z^\theta$ is a scheme of finite type (which will be proven
shortly) the intersection $\Gr_P^\theta\cap \Gr_{U(P^-)}$ is also a scheme of finite type.

\section{Relation of the Zastava spaces with $\BunPbw$}  \label{relation with BunPbw}

Our goal now is to show that the space $Z^\theta$ models the stack $\BunPbw$ from the point of view of singularities.

\ssec{}

First, we have to introduce the following relative version of $Z^{\theta}$. 

Let $\F^b_M$ be a fixed $M$-bundle on $X$ and let $\F^b_G$ be the induced 
$G$-bundle under our fixed embedding $M\hookrightarrow G$.
The space $Z^\theta_{\F^b_M}$ is defined as follows: it classifies quadruples 
$(\F_G,\F_M,\beta_M,\beta)$ as in the case of $Z^\theta$
with the only difference that the trivial $M$-bundle $\F^0_M$ is replaced by 
$\F^b_M$ and the trivial $G$-bundle $\F^0_G$
is replaced by $\F^b_G$. 

Since every $M$-bundle is locally trivial (cf. \cite{DS}),
and due to the Factorization property, the spaces
$Z^\theta$ and $Z^\theta_{\F^b_M}$ are \'etale-locally isomorphic.

Similarly, if $S$ is a scheme (or a stack) mapping to $\Bun_M$, we can 
define the space $Z^\theta_S$. When $S$ is smooth, then using \cite{DS} we
obtain that $Z^\theta_S$ is locally in the smooth topology equivalent to $Z^\theta$. 
In practice, we will take $S$ to be $\Bun_M$.

Along the same lines, we define the relative version 
$\Mod_{M,S}^{+,\theta}$ of $\Mod_M^{+,\theta}$. 
When $S=\Bun_M$, we will denote it by $\Mod_{\Bun_M}^{+,\theta}$.

Recall that the stack $\Bun_M$ splits into connected 
components numbered by the elements of $\Lambda_{G,P}$. By definition,
a point $\F_M$ belongs to the connected component $\Bun_M^\theta$ if the associated $M/[M,M]$-bundle is of degree $-\theta$. 
We will use the superscript $\theta$ to designate the corresponding connected component of the stack $\BunPbw$ or
$\Bun_{P^-}$. 

\begin{prop} \label{identification}
For every $\theta\in\Lambda_{G,P}^{\on{pos}}$ and $\theta'\in\Lambda_{G,P}$
there is a canonical isomorphism between $Z^{\theta}_{\Bun_M}$ and an open sub-stack of
$\BunPbw^{\theta+\theta'}\underset{\Bun_G}\times \Bun^{\theta'}_{P^-}$. 
\end{prop}

From this proposition it follows, in particular, that $Z^{\theta}_{\Bun_M}$ is a stack of finite type,
and hence $Z^\theta$ is a scheme of finite type (and not just an ind-scheme). 

\begin{proof}

Let us analyze what it means to have an $S$-point of the Cartesian product 
$$\BunPbw^{\theta+\theta'}\underset{\Bun_G}\times \Bun^{\theta'}_{P^-}.$$
By definition, we have a $G$-bundle, a pair of $M$-bundles $\F_M$ and $\F'_M$
and two systems of maps $\kappa$ and $\kappa^-$ for every $G$-module $\V$:
\begin{align*}
&\kappa:(\V^{U(P)})_{\F_M} \to \V_{\F_G}  \\
&\kappa^-:\V_{\F_G} \to (\V_{U(P^-)})_{\F'_M},
\end{align*}
which satisfy the Pl\"ucker relations, with the condition that the $\kappa^-$'s are surjective,
and the $\kappa$'s are injective over every geometric point of $S$.

\medskip

We define the open sub-stack $(\BunPbw^{\theta+\theta'}\underset{\Bun_G}\times \Bun^{\theta'}_{P^-})^0$ by
the following condition: for every geometric point $s\in S$, the $P$- and $P^-$-structures defined on $\F_G|_s$
by means of $\kappa$ and $\kappa^-$ {\it at the generic point of $X$} are mutually transversal.

Let $(\F_G,\F_M,\F'_M,\beta_M,\beta)$ be an $S$-point of $Z^\theta_{\Bun_M}$. It is clear that the maps
$$(\V^{U(P)})_{\F_M}\to \V_{\F_G} \text{ and }\V_{\F_G}\to (\V_{U(P^-)})_{\F'_M},$$ 
as in the definition of $Z^\theta_{\Bun_M}$
indeed define an $S$-point of $(\BunPbw^{\theta+\theta'}\underset{\Bun_G}\times \Bun^{\theta'}_{P^-})^0$.

Conversely, given an $S$-point of $(\BunPbw^{\theta+\theta'}\underset{\Bun_G}\times \Bun^{\theta'}_{P^-})^0$,
we define an $S$-point of $Z^\theta_{\Bun_M}$ as follows: 

Firstly, we set the ``background'' $M$-bundle to be $\F'_M$. Let $\F'_G$ denote the induced $G$-bundle under
$M\hookrightarrow G$. Secondly, by construction, 
there exists an open dense subset $(X\times S)^0\subset X\times S$, such that $\F_G|_{(X\times S)^0}$ admits reductions
simultaneously to $P$ and $P^-$, which are, moreover, transversal. Hence, over $(X\times S)^0$, we have identifications
$\beta:\F_G\simeq \F'_G$ and $\beta_M:\F_M\simeq \F'_M$. Therefore, it remains to show that $\beta_M$ is such that the maps
$\beta_M^\U:\U_{\F_M}\to \U_{\F'_M}$, which are defined on $(X\times S)^0$, 
extend as regular maps to the entire $X\times S$,
provided that $\U$ is of the form $\V^{U(P)}$ for a $G$-module $\V$.

However, note that for $\U$ of the above form the composition 
$$\U\to \on{Ind}(\U)^{U(P)}\to \on{Ind}(\U)\to \on{Ind}(\U)_{U(P^-)}$$
is an isomorphism. (As a remark, let us observe that in char 0, the maps $\V^{U(P)}\to \V_{U(P^-)}$ 
are isomorphisms for any $G$-module $\V$, which is not necessarily the case over $\Fq$.)
Hence, by composing $\kappa$ and $\kappa^-$ we do obtain a regular map
$$\U_{\F_M}\to (\on{Ind}(\U)^{U(P)})_{\F_M}\to \V_{\F_G}\to (\on{Ind}(\U)_{U(P^-)})_{\F'_M}\simeq \U_{\F'_M},$$
which is what we had to show.

\end{proof}

\ssec{}

Observe that under the isomorphism of the above proposition, the open sub-stack ${Z_{max}^\theta}_{\Bun_M}$ coincides with 
the preimage of $\Bun_P\subset \BunPbw$. Let us now analyze the behavior of other strata of $\BunPbw$ under the isomorphism
of \propref{identification}.

Recall that for $\theta\in\Lambda_{G,P}^{\on{pos}}$ we introduced a locally 
closed sub-stack $_\theta\BunPb\subset \BunPb$, as the image of the natural map
$$X^\theta\times \Bun_P\to\BunPb.$$

\medskip

Let $_\theta\BunPbw\subset \BunPbw$ denote the preimage of $_\theta\BunPb$ under $\r:\BunPbw\to\BunPb$.
As in \propref{descrBunPbw} one shows that
$$_\theta\BunPbw\simeq \Bun_P\underset{\Bun_M}\times \Mod_{\Bun_M}^{+,\theta}.$$

\bigskip

Let $(_{\theta'}\BunPbw^{\theta+\theta'}\underset{\Bun_G}\times \Bun^{\theta'}_{P^-})^0$ be the preimage
of $_{\theta'}\BunPbw^{\theta+\theta'}$ under the natural projection. 

\begin{lem}  \label{identplus}
The stack 
$(_{\theta'}\BunPbw^{\theta+\theta'}\underset{\Bun_G}\times\Bun^{\theta'}_{P^-})^0$ is empty unless 
$\theta-\theta'\in \Lambda_{G,P}^{\on{pos}}$. For $\theta'=\theta$, the above sub-stack identifies with the image
of $\Mod_{\Bun_M}^{+,\theta}$ under $\s^\theta:\Mod_{\Bun_M}^{+,\theta}\hookrightarrow Z^\theta_{\Bun_M}$.
\end{lem}

\begin{proof}

Note that an $S$-point of $\BunPb$ belongs to $_{\theta'}\BunPb$ if and only if the following condition holds:
for every $G$-dominant weight $\lambdach$, orthogonal to $\on{Span}(\alpha_\i)$ for $\i\in \I_M$, the corresponding map
$$\kappa^{\lambdach}:\L^{\lambdach}_{\F_{M/[M,M]}}\to \V^{\lambdach}_{\F_G}$$
is such that there exists a short exact sequence
$$0\to \M_1\to \on{coker}(\kappa^{\lambdach})\to\M_2\to 0,$$
such that $\M_2$ is a vector bundle on $X\times S$, the support of $\M_1$ is $X$-finite and over 
any geometric point $s\in S$, the length of $\M_1|_s$ is exactly $\langle \theta',\lambdach\rangle$.

Given an $S$-point of $(_{\theta'}\BunPbw^{\theta+\theta'}\underset{\Bun_G}\times\Bun^{\theta'}_{P^-})^0$,
we can compose the above embedding of sheaves with
$$\V^{\lambdach}_{\F_G}\overset{\kappa^-}\longrightarrow (\V^{\lambdach}_{U(P^-)})_{\F'_M}
\to \L^{\lambdach}_{\F'_{M/[M,M]}}$$ 
and we obtain a map between line bundles, such that over every
geometric point $s\in S$ its total amount of zeroes is $\langle\theta,\lambdach\rangle$. This readily
implies the first assertion of the lemma.

\medskip

To prove the second assertion, observe that
$$(_\theta\BunPbw^{\theta+\theta'}\underset{\Bun_G}\times\Bun^{\theta'}_{P^-})^0\simeq 
(\Bun^{\theta'}_{P^-}\underset{\Bun_G}\times\Bun^{\theta'}_P)^0\underset{\Bun_M}\times \Mod_{\Bun_M}^{+,\theta}.$$
However, the condition on the degree forces that 
$$(\Bun^{\theta'}_{P^-}\underset{\Bun_G}\times\Bun^{\theta'}_P)^0\simeq \Bun_M.$$

Hence,
$$(_\theta\BunPbw^{\theta+\theta'}\underset{\Bun_G}\times\Bun^{\theta'}_{P^-})^0\simeq \Mod_{\Bun_M}^{+,\theta}$$
and the fact that its embedding into $Z^\theta_{\Bun_M}$ coincides with $\s^\theta$ follows from the construction.

\end{proof}

\ssec{}   \label{various strata}

For an element $\theta'$ with $\theta-\theta'\in \Lambda_{G,P}^{\on{pos}}$, let 
us denote by $_{\theta'}Z^\theta$
the corresponding locally closed subvariety of $Z^\theta$, which is the trace of $_{\theta'}\BunPbw$
under the isomorphism of \propref{identification}. In particular, $_{0}Z^\theta=Z^\theta_{max}$.

\medskip

As in \propref{descrBunPbw} above, we obtain:
$$_{\theta'}Z^\theta\simeq Z_{max}^{\theta-\theta'}\underset{\Bun_M}\times \Mod_M^{+,\theta'},$$
where the map $Z_{max}^{\theta-\theta'}\to \Bun_M$ used in the definition of the fiber product
is $(\F_M,\beta_M)\in Z^{\theta-\theta'}\mapsto \F_M\in \Bun_M$.

\bigskip

Let us denote by $_{\theta'}\SS^\theta$ 
the intersection of $_{\theta'}Z^\theta$ with
the central fiber $\SS^\theta$. 
We obtain the following description of $_{\theta'}\SS^\theta$:

\medskip

Let $\Conv_M$ denote the convolution diagram  of the affine Grassmannian of the group $M$. By definition,
$\Conv_M$ classifies quadruples $(\F_M,\F'_M,\wt\beta_M,\beta'_M)$, where $\wt\beta_M$ is an isomorphism 
$\F_M|_{X-x}\simeq \F'_M|_{X-x}$ and $\beta'_M$ is an isomorphism $\F'_M|_{X-x}\simeq \F^0_M|_{X-x}$. We have a natural
projection $pr':\Conv_M\to\Gr_M$, which sends $(\F_M,\F'_M,\wt\beta_M,\beta'_M)\mapsto (\F'_M,\beta'_M)$ and the projection
$pr:\Conv_M\to\Gr_M$, which sends $(\F_M,\F'_M,\wt\beta_M,\beta'_M)\mapsto (\F_M,\beta'_M\circ \wt\beta_M)$.

Projection $pr'$ makes $\Conv_M$ a fibration over $\Gr_M$ with the typical fiber isomorphic to $\Gr_M$ and
we will denote by $\Conv_M^{+,\theta}$ the closed sub-scheme of $\Conv$, which is a fibration
over $\Gr_M$ with the typical fiber $\Gr_M^{+,\theta}$.

Using \propref{central fiber} we obtain:
$$_{\theta'}\SS^\theta\simeq (\Gr_P^{\theta-\theta'}\cap \Gr_{U(P^-)})\underset{\Gr_M}\times \Conv_M^{+,\theta'},$$
where $\Gr_M\leftarrow \Conv_M^{+,\theta'}$ is the map $pr'$.

\ssec{Smoothness issues}  \label{smoothness}

Above we have constructed the map
$$Z^{\theta}_{\Bun_M}\simeq (\BunPbw^{\theta+\theta'}\underset{\Bun_G}\times\Bun^{\theta'}_{P^-})^0\to 
\BunPbw^{\theta+\theta'}.$$
We do not {\it a priori} know whether this map is smooth, since $\Bun_{P^-}\to\Bun_G$ is not smooth. We will now construct
an open sub-stack in $Z^{\theta}_{\Bun_M}$, which will map smoothly onto $\BunPbw^{\theta+\theta'}$.

\medskip

Let ${\mathfrak u}(P)$ be the Lie algebra of $U(P)$ viewed as an $M$-module. We define the open sub-stack $\Bun_M^r\subset \Bun_M$
to consist of those $M$-bundles $\F_M$, for which $H^1(X,\U)=0$, for all $M$-modules
$\U$, which appear as sub-quotients of ${\mathfrak u}(P)$.
Let $\Bun_{P^-}^r$ be the preimage of $\Bun_M^r$ under the natural projection $\qo:\Bun_{P^-}\to \Bun_M$. 

\begin{lem}  
The restriction of the natural map $\Bun_{P^-}\to\Bun_G$ to $\Bun_{P^-}^r$ is smooth.
\end{lem}

\begin{proof}

Since both $\Bun_{P^-}$ and $\Bun_G$ are smooth, it is enough to check the surjectivity
on the level of tangent spaces. Thus, let $\F_{P^-}$ be a $P^-$-bundle and let $\F_G$ be the induced $G$-bundle.
We must show that
$$H^1(X,{\mathfrak p}^-_{\F_{P^-}})\to H^1(X,{\mathfrak g}_{\F_G})$$
is surjective if $\qo(\F_{P^-})\in \Bun_M^r$.

In general, the cokernel of this map is $H^1(X,({\mathfrak g}/{\mathfrak p}^-)_{\F_{P^-}})$. However, 
the irreducible subquotients of
${\mathfrak g}/{\mathfrak p}^-$ as a $P^-$-module are all $M$-modules, which appear 
in the Jordan-Holder series of ${\mathfrak u}(P)$.
Hence, the assertion of the proposition follows from the definition of $\Bun_M^r$.

\end{proof}

Let $Z^\theta_{\Bun_M^r}$ denote the corresponding open sub-stack of $Z^\theta_{\Bun_M}$. From \propref{identification}
we obtain an isomorphism $Z^\theta_{\Bun_M^r}\simeq (\BunPbw^{\theta+\theta'}\underset{\Bun_G}\times \Bun^{\theta',r}_{P^-})^0$
and from the above lemma, we see that the resulting map $Z^\theta_{\Bun_M^r}\to \BunPbw^{\theta+\theta'}$ is smooth.
In particular, since the stack $\Bun_P$ is smooth, we obtain the following corollary:

\begin{cor}  \label{opensmooth}
The open sub-scheme $Z^\theta_{max}$ of $Z^\theta$ is smooth.
\end{cor}

It is well-known (cf. \cite{DS}) that every open sub-stack of $\Bun_G$ of finite type belongs to the 
image of some $\Bun^{\theta',r}_{P^-}$, when $-\theta'$ is large enough. Similarly, it is easy to see that
every open sub-stack of finite type of $\BunPbw^{\theta}$ eventually belongs to the image of $Z^{\theta-\theta'}_{\Bun_M^r}$.
Hence, in order to understand the singularities of $\BunPbw$, it is sufficient to analyze the singularities of $Z^{\theta}$.

\section{Computation of $\IC_{Z^\theta}$: Statements}

\ssec{}

For $\theta\in \Lambda^{\on{pos}}_{G,P}\simeq \on{Span}^+(\alpha_\i,\,\i\in \I-\I_M)$, let
${\mathfrak P}(\theta)$ denote an element of the set of partitions of $\theta$ as a sum
$\theta=\underset{k}\Sigma\, \theta_k$, where each $\theta_k$ is a projection under $\Lambda\to\Lambda_{G,P}$ of a
coroot of $G$ belonging to $\on{Span}^+(\alpha_\i,\,\i\in \I-\I_M)$.

We emphasize the difference between ${\mathfrak P}(\theta)$ and ${\mathfrak A}(\theta)$: in the latter case we
decompose $\theta$ as a sum of arbitrary non-zero elements of $\Lambda^{\on{pos}}_{G,P}$.

For a fixed ${\mathfrak P}(\theta)$, let $X^{{\mathfrak P}(\theta)}$ denote the corresponding partially symmetrized
power of the curve. In other words, if $\theta=\underset{k}\Sum\, n_k\cdot \theta_k$, where $\theta_k$'s are pairwise
distinct, $X^{{\mathfrak P}(\theta)}=\underset{k}\Pi\, X^{(n_k)}$.

Now we need to introduce a version of the Beilinson-Drinfeld affine Grassmannian $\Gr_M^{{\mathfrak P}(\theta)}$.
First, consider the ind-scheme $\Gr_M^{{\mathfrak P}(\theta),\infty}$, which classifies triples 
$(x^{{\mathfrak P}(\theta)},\F_M,\beta_M)$, 
where $x^{{\mathfrak P}(\theta)}\in X^{{\mathfrak P}(\theta)}$,
$\F_M$ is an $M$-bundle on $X$ and $\beta_M$ is the trivialization of $\F_M$ away from the support of $x^{{\mathfrak P}(\theta)}$.
(We leave it to the reader to formulate the above definition in terms of $S$-points, in the spirit of what we have
done before.)

Consider the open subset $\overset{o}X{}^{{\mathfrak P}(\theta)}$ of $X^{{\mathfrak P}(\theta)}$ equal to the complement 
of all the diagonals. Inside $\Gr_M^{{\mathfrak P}(\theta),\infty}|_{\overset{o}X{}^{{\mathfrak P}(\theta)}}$ we define
the closed subset $\Gr_M^{{\mathfrak P}(\theta)}|_{\overset{o}X{}^{{\mathfrak P}(\theta)}}$ as follows:
For $x^{{\mathfrak P}(\theta)}=\Sum \, \theta_k\cdot x_k$ with all the $x_k$'s distinct, the fiber of $\Gr_M^{{\mathfrak P}(\theta),\infty}$
over it is just the product of the affine Grassmannians $\underset{k}\Pi\, \Gr_{M,x_k}$ and the fiber of
$\Gr_M^{{\mathfrak P}(\theta)}$ is set to be $\underset{k}\Pi\, \grb^{\flat(\theta_k)}_{M,x_k}$, where $\flat(\theta_k)$ is
as in \propref{pospart}. The entire $\Gr_M^{{\mathfrak P}(\theta)}$ is defined as a closure of 
$\Gr_M^{{\mathfrak P}(\theta)}|_{\overset{o}X{}^{{\mathfrak P}(\theta)}}$ inside $\Gr_M^{{\mathfrak P}(\theta),\infty}$.

\medskip

By construction, if $(x^{{\mathfrak P}(\theta)},\F_M,\beta_M)$ belongs to $\Gr_M^{{\mathfrak P}(\theta)}$, then among the rest, the
trivialization $\beta_M$ has the following property: for every $G$-module $\V$ the map
$$\beta_M^{\V^{U(P)}}:(\V^{U(P)})_{\F_M}\to (\V^{U(P)})_{\F^0_M}\simeq \V^{U(P)}\otimes \O_X,$$
which is defined {\it a priori} on $X-x^{{\mathfrak P}(\theta)}$ extends to a regular map on $X$. Therefore, we obtain a map
$i_{{\mathfrak P}(\theta)}:\Gr_M^{{\mathfrak P}(\theta)}\to \Mod_M^{+,\theta}$, which covers the natural map
$X^{{\mathfrak P}(\theta)}\to X^\theta$. It is easy to see that the above map $i_{{\mathfrak P}(\theta)}$ is finite.

\ssec{}

Let us denote by $\IC^{{\mathfrak P}(\theta)}$ the intersection cohomology sheaf on $\Gr_M^{{\mathfrak P}(\theta)}$.
We need to understand more explicitly the behavior of  $\IC^{{\mathfrak P}(\theta)}$ over the diagonals in 
$X^{{\mathfrak P}(\theta)}$.

Thus, let $\Delta_X\subset X^{{\mathfrak P}(\theta)}$ be the main diagonal. By construction, 
$\Gr_M^{{\mathfrak P}(\theta)}|_{\Delta_X}$ is a sub-scheme of the relative affine Grassmannian $\Gr_{M,X}$. Recall that
$\on{Loc}_X$ denotes the localization functor from $\on{Rep}(\check M)$ to the category of perverse sheaves on $\Gr_{M,X}$.

\begin{lem}  \label{restrM}
If ${\mathfrak P}(\theta)$ corresponds to $\theta=\underset{k}\Sum\, n_k\cdot \theta_k$, then
the *-restriction of $\IC^{{\mathfrak P}(\theta)}$ to $\Gr_M^{{\mathfrak P}(\theta)}|_{\Delta_X}\subset \Gr_{M,X}$
can be canonically identified with 
$$\on{Loc}_X(\underset{k}\otimes\, \on{Sym}^{n_k}(\check{\mathfrak u}(P)_{\theta_k}))\otimes 
(\Ql(\frac{1}{2})[1])^{\otimes |{\mathfrak P}(\theta)|-1},$$
where $|{\mathfrak P}(\theta)|=\underset{k}\Sum\, n_k$.
\end{lem}

\begin{proof}

Consider the corresponding non-symmetrized power of the curve $\underset{k}\Pi\, X^{n_k}$. Over it we can consider
the scheme $\underset{k}\Pi\, (\grb_{M,X}^{\flat(\theta_k)})^{n_k}$ and we have a natural proper map
$$sym:\underset{k}\Pi\, (\grb_{M,X}^{\flat(\theta_k)})^{n_k}\to \Gr_M^{{\mathfrak P}(\theta)},$$
which covers the usual symmetrization map $\underset{k}\Pi\, X^{n_k}\to X^{{\mathfrak P}(\theta)}$.
Let us denote temporarily by $\S$ the direct image $sym_!(\IC_{\underset{k}\Pi\, (\grb_{M,X}^{\flat(\theta_k)})^{n_k}})$.

The fact that the map which defines
convolution of perverse sheaves on the usual affine Grassmannian $\Gr_M$ is semi-small implies that the above map
$sym$ is small. Hence, $\S$ is the Goresky-MacPherson extension of its restriction to the open sub-scheme 
$\Gr_M^{{\mathfrak P}(\theta)}|_{\overset{o}X{}^{{\mathfrak P}(\theta)}}$. In particular, it carries a canonical
action of the product of symmetric groups $\underset{k}\Pi\, S^{n_k}$, because this is obviously so over 
$\Gr_M^{{\mathfrak P}(\theta)}|_{\overset{o}X{}^{{\mathfrak P}(\theta)}}$, and $\IC^{{\mathfrak P}(\theta)}\simeq (\S)^{\underset{k}\Pi\,
S^{n_k}}$.

By construction, the *-restriction of $\S$ to $\Gr_M^{{\mathfrak P}(\theta)}|_{\Delta_X}$ can be identified with
$$\on{Loc}_X(\underset{k}\otimes\, (\check{\mathfrak u}(P)_{\theta_k})^{\otimes n_k})\otimes
(\Ql(\frac{1}{2})[1])^{\otimes \underset{k}\Sigma\, n_k-1}.$$

Therefore, it remains to see that the $\underset{k}\Pi\, S^{n_k}$-action on $\S|_{\Delta_X}$ corresponds to the natural
action of the group on $\underset{k}\otimes\,(\check{\mathfrak u}(P)_{\theta_k})^{\otimes n_k}$. We prove the
latter fact as follows:

Since taking the global cohomology is a fiber functor for the category of spherical perverse sheaves on $\Gr_M$, it suffices to 
analyze the $\underset{k}\Pi\, S^{n_k}$-action on the direct image of $\S$ under 
$\Gr_M^{{\mathfrak P}(\theta)}\to X^{{\mathfrak P}(\theta)}$,
in which case the assertion becomes obvious.

\end{proof}

\ssec{The main theorem}

Our main technical result is the following theorem:

\begin{thm}  \label{calcmain}
The $!$-restriction of $\IC_{Z^{\theta}}$ under $\s^\theta:\Mod_M^{+,\theta}\to Z^\theta$ can be identified with
$$\underset{{\mathfrak P}(\theta)}\oplus\, i_{{\mathfrak P}(\theta)}{}_*(\IC^{{\mathfrak P}(\theta)})\otimes
(\Ql(\frac{1}{2})[1])^{-|{\mathfrak P}(\theta)|}.$$
\end{thm}

\noindent{\it Remark.}
Let us explain to what extent the isomorphism stated in this theorem is canonical. (In fact, it is not!)
The LHS carries the cohomological filtration ({\it filtration canonique}), which corresponds
to the filtration on the RHS according to $|{\mathfrak P}(\theta)|$
under whichever isomorphism between the RHS and the LHS that we choose. Unfortunately, our proof does not 
give a canonical identification even for the associated graded quotients: each
$i_{{\mathfrak P}(\theta)}{}_*(\IC^{{\mathfrak P}(\theta)})$ appears up to tensoring with a $1$-dimensional
vector space.

\medskip

To prove this theorem, we will proceed by induction on $|\theta|:=\underset{\i\in \I-\I_M}\Sum\, n_\i$ if 
$\theta=\underset{\i\in \I-\I_M}\Sum\, n_\i\cdot \alpha_\i$. First, we will derive from it various facts 
about $\IC_{\BunPbw}$, since we will use them to perform the induction step.

\ssec{}

Observe that since $\IC_{Z^\theta}$ is Verdier self-dual, from \thmref{calcmain} we obtain the description of
$\s^{\theta}{}^*(\IC_{Z^{\theta}})$ as well. By translating this description to $\BunPbw$ using \propref{identification},
we obtain the following corollary:

\begin{cor}  \label{restr1}
The *-restriction of $\IC_{\BunPbw}$ to $_\theta\BunPbw\simeq \Bun_P\underset{\Bun_M}\times \Mod_{\Bun_M}^{+,\theta}$ 
can be identified with
$$\underset{{\mathfrak P}(\theta)}\oplus\, (\id\times i_{{\mathfrak P}(\theta)})_*
(\IC_{\Bun_P\underset{\Bun_M}\times \H_M^{{\mathfrak P}(\theta)}})\otimes
(\Ql(\frac{1}{2})[1])^{|{\mathfrak P}(\theta)|},$$
where $\H_M^{{\mathfrak P}(\theta)}$ is the corresponding relative version of $\Gr_M^{{\mathfrak P}(\theta)}$ over $\Bun_M$.
\end{cor}

From this corollary one easily deduces \thmref{restr}:

\begin{proof} (of \thmref{restr})

To simplify the notation, we will take the element 
${\mathfrak A}(\theta)={\mathfrak A}^0(\theta)$ corresponding to the decomposition
which consists of one element: $\theta=\theta$. We need to calculate the *-restriction of $\IC_{\BunPbw}$
to $$_{{\mathfrak A}^0(\theta)}\BunPbw\simeq
\Bun_P\underset{\Bun_M}\times \H_M^{{\mathfrak A}^0(\theta)}\simeq 
\Bun_P\underset{\Bun_M}\times \H_{M,X}^{+,\theta}.$$

By definition, our embedding $j_{{\mathfrak A}^0(\theta)}$ factors through 
$\Bun_P\underset{\Bun_M}\times \Mod_{\Bun_M}^{+,\theta}\simeq{}_\theta\BunPbw$.
Note that $\H_{M,X}^{+,\theta}\subset \Mod_{\Bun_M}^{+,\theta}$ is exactly the preimage of the main diagonal 
$\Delta_X\subset X^\theta$.
Therefore, the sought-for complex is, according to \corref{restr1}, the direct sum over ${\mathfrak P}(\theta)$ of
$$(\id\times i_{{\mathfrak P}(\theta)})_*
(\IC_{\Bun_P\underset{\Bun_M}\times \H_M^{{\mathfrak P}(\theta)}})|_{\Delta_X}
\otimes(\Ql(\frac{1}{2})[1])^{|{\mathfrak P}(\theta)|}.$$

Using \lemref{restrM}, we obtain that 
$(\id\times i_{{\mathfrak P}(\theta)})_*(\IC_{\Bun_P\underset{\Bun_M}
\times \H_M^{{\mathfrak P}(\theta)}})|_{\Delta_X}$
corresponding to ${\mathfrak P}(\theta)$ with $\theta=\underset{k}\Sum\, n_k\cdot \theta_k$ equals 
$$\on{Loc}_{\Bun_P,X}(\underset{k}\otimes\, \on{Sym}^{n_k}(\check{\mathfrak u}(P)_{\theta_k}))\otimes 
(\Ql(1)[2])^{\otimes |{\mathfrak P}(\theta)|})\otimes (\Ql(\frac{1}{2})[1])^{\otimes -1}.$$
However,
$$\underset{{\mathfrak P}(\theta)}\oplus (\underset{k}\otimes\, \on{Sym}^{n_k}(\check{\mathfrak u}(P)_{\theta_k})
\otimes (\Ql(1)[2])^{\otimes |{\mathfrak P}(\theta)|})\simeq 
\underset{i\geq 0}\oplus\,\on{Sym}^i(\check{\mathfrak u}(P))_\theta\otimes(\Ql(1)[2])^{\otimes i} ,$$
which is what we had to show.

\end{proof}

\ssec{}

The following result is an interesting byproduct of \corref{restr1}. In order to save notation, we will formulate it
for ${\mathfrak A}(\theta)={\mathfrak A}^0(\theta)$, although the generalization to an arbitrary ${\mathfrak A}(\theta)$
is straightforward.

Consider the {\it hyperbolic} restriction of $\IC_{\BunPbw}$ to a stratum $_{{\mathfrak A}^0(\theta)}\BunPbw$. By definition,
this is the !-restriction of $\IC_{\BunPbw}$ to $_\theta\BunPbw$ followed by the further *-restriction from 
$_\theta\BunPbw$ to $_{{\mathfrak A}^0(\theta)}\BunPbw$. 

\begin{cor}  \label{restr2}
The hyperbolic restriction (in the above sense) of $\IC_{\BunPbw}$ to \newline
$_{{\mathfrak A}^0(\theta)}\BunPbw\simeq \Bun_P\underset{\Bun_M}\times \H_{M,X}^{+,\theta}$ is isomorphic to
$$\on{Loc}_{\Bun_P,X}(\on{Sym}(\check{\mathfrak u}(P))_{\theta})\otimes 
(\Ql(\frac{1}{2})[1])^{\otimes -1}.$$
\end{cor}

Let us draw the reader's attention to the fact that \corref{restr2} implies that the hyperbolic restriction
of $\IC_{\BunPbw}$ to $_{{\mathfrak A}^0(\theta)}\BunPbw$ is a perverse sheaf, up to a cohomological shift.

The proof of this corollary repeats the above proof of \thmref{restr}, using the fact that 
$$\underset{{\mathfrak P}(\theta)=\underset{k}\sum\, n_k\theta_k} \oplus\, \underset{k}\otimes\,
\on{Sym}^{n_k}(\check{\mathfrak u}(P)_{\theta_k})\simeq
\on{Sym}(\check{\mathfrak u}(P))_{\theta}.$$

\noindent{\it Remark.}
We remark again that, due to the non-canonicity of the direct sum decomposition
stated in \thmref{calcmain}, the isomorphism of \corref{restr2} is non canonical either. We only
can claim that the LHS carries a canonical filtration, which on the RHS coincides with the filtration by the degree.
However, in the course of the proof of \thmref{calcmain}, we will show that the above hyperbolic restriction can
be {\it canonically} identified with $\on{Loc}_X(U(\check{\mathfrak u}(P))_{\theta})\otimes 
(\Ql(\frac{1}{2})[1])^{\otimes -1}$. It seems natural to guess (although our proof does not imply it)
that our filtration on the LHS corresponds under this isomorphism to the canonical filtration on $U(\check{\mathfrak u}(P))$.

\section{Computation of $\IC_{Z^\theta}$: Proofs}

\ssec{}

The goal of this section is to prove \thmref{calcmain}. Our strategy will be as follows: from the induction
hypothesis we will obtain an almost complete description of how $\s^{\theta}{}^!(\IC_{Z^\theta})$ looks
like away from the main diagonal. Then we will explicitly compute the 
``contribution" at the main diagonal and hence
prove the theorem.

The crucial idea of the proof is the following assertion:

\begin{prop}  \label{directandinverse}
There is a canonical isomorphism $\s^{\theta}{}^!(\IC_{Z^\theta})\simeq \pi_P{}_!(\IC_{Z^\theta})$.
\end{prop}

We will deduce this proposition from the following well-known lemma:

Let $\pi:Y'\to Y$ be a map of schemes and $\s:Y\to Y'$ be a section. 
Assume now that the group ${\mathbb G}_m$ acts on
$Y'$ in such a way that it 
``contracts" $Y'$ onto $Y$. This means that the action map ${\mathbb G}_m\times Y'\to Y'$
extends to a regular map ${\mathbb A}^1\times Y'\to Y'$, such that the composition
$$0\times Y'\to {\mathbb A}^1\times Y'\to Y'$$ coincides with $\s\circ \pi:Y'\to Y\to Y'$. Let now $\S$ be a 
${\mathbb G}_m$-equivariant complex on $Y'$.

\begin{lem}
Under the above circumstances, $\pi_!(\S)\simeq \s^!(\S)$.
\end{lem}

We apply this lemma to $Y'=Z^\theta$ and $Y=\Mod_M^{+,\theta}$. To construct a ${\mathbb G}_m$-action we proceed as follows.
Let ${\mathbb G}_m\to Z(M)$ be a 1-parameter subgroup, which acts as contraction on $U(P^-)$. In this way
${\mathbb G}_m$ acts on the trivial $M$-bundle $\F^0_M$ on $X$ and hence on $Z^\theta$. It remains to verify that
${\mathbb G}_m$ indeed contracts $Z^\theta$ onto $\Mod_M^{+,\theta}$. We do that as follows:

\medskip

Let $\Gr_G^{\theta,\infty}$ be the Beilinson-Drinfeld affine Grassmannian over $X^\theta$. In other words, a point of
$\Gr_G^{\theta,\infty}$ is a triple $(x^\theta,\F_G,\beta)$, where $x^\theta\in X^\theta$, $\F_G$ is a $G$-bundle on $X$
and $\beta$ is a trivialization of $\F_G$ off the support of $x^\theta$. We have the unit section $X^\theta\to \Gr_G^{\theta,\infty}$
which sends $x^\theta$ to $(x^\theta,\F^0_G,\beta^0)$, where $\beta^0$ is the tautological trivialization of the trivial bundle.

In the same way we can consider a Beilinson-Drinfeld version of the affine Grassmannian for the group $U(P^-)$ (denote it
$\Gr_{U(P^-)}^{\theta,\infty}$), and we have a locally closed embedding: 
$\Gr_{U(P^-)}^{\theta,\infty}\hookrightarrow \Gr_G^{\theta,\infty}$. 

By construction, our $Z^\theta$ is a closed sub-scheme
inside $\Mod_M^{+,\theta}\underset{X^\theta}\times \Gr_{U(P^-)}^{\theta,\infty}$. The image of $\s^\theta$ is
the product of $\Mod_M^{+,\theta}$ and the unit section of $\Gr_{U(P^-)}^{\theta,\infty}$.

The above ${\mathbb G}_m$-action on $Z^\theta$ comes from a natural action of this group on $\Gr_{U(P^-)}^{\theta,\infty}$, 
while the action on $\Mod_M^{+,\theta}$ is trivial. Therefore, to prove our assertion we have to show that 
${\mathbb G}_m$ contracts $\Gr_{U(P^-)}^{\theta,\infty}$ to the unit section. However, this easily follows from the fact
that our ${\mathbb G}_m\to Z(M)$ contracts $U(P^-)$ to $1\in U(P^-)$.

\ssec{}

Having established \propref{directandinverse}, we obtain the following corollary:

\begin{cor}    \label{semisimple}
When we pass from $\Fq$ to $\Fqb$, the complex $\s^{\theta}{}^!(\IC_{Z^\theta})\simeq \pi_P{}_!(\IC_{Z^\theta})$
splits as a direct sum of (cohomologically shifted) irreducible perverse sheaves.
\end{cor}

\begin{proof}

According to \cite{BBD}, $\pi_P{}_!(\IC_{Z^\theta})$ has weights $\leq 0$, since $\IC_{Z^\theta}$
is pure of weight $0$. At the same time, \cite{BBD} implies that $\s^{\theta}{}^!(\IC_{Z^\theta})$ has weights
$\geq 0$. Hence, we obtain 
that $\s^{\theta}{}^!(\IC_{Z^\theta})\simeq \pi_P{}_!(\IC_{Z^\theta})$
is pure of weight $0$.

Hence, the assertion of the corollary follows from the decomposition theorem. 

\end{proof}

From this moment until \secref{Weil}, we will disregard the $\Fq$-structure on $\s^{\theta}{}^!(\IC_{Z^\theta})$ and
will prove the isomorphism stated in \thmref{calcmain} over $\Fqb$. In \secref{Weil} we will show that the direct
sum decomposition holds over $\Fq$ as well.

\medskip

Now let us use the induction hypothesis, i.e. our knowledge about $\s^{\theta'}{}^!(\IC_{Z^{\theta'}})$
for all $\theta'$ with $\theta-\theta'\in\Lambda_{G,P}^{\on{pos}}$. It is easy to see that the {\it Factorization property} 
of \propref{factor} implies that locally over $X^\theta-\Delta_X$ we do obtain an isomorphism
$$\s^{\theta}{}^!(\IC_{Z^\theta})\simeq \underset{{\mathfrak P}(\theta),\,|{\mathfrak P}(\theta)|\neq 1}\oplus\, 
i_{{\mathfrak P}(\theta)}{}_*(\IC^{{\mathfrak P}(\theta)})[-|{\mathfrak P}(\theta)|].$$

However, globally we can {\it a priori} have a non-trivial monodromy: suppose that ${\mathfrak P}(\theta)$ corresponds
to $\theta=2\cdot \theta'$, where $\theta'$ is the image of a coroot in $\on{Span}(\alpha_j,\,\j\in\I-\I_M)$.
Then the preimage of $X^\theta-\Delta_X$ in $X^{{\mathfrak P}(\theta)}$ 
is $X^{(2)}-\Delta_X$, and we can have an order $2$ monodromy. Therefore, so far we can only claim that 
$$\s^{\theta}{}^!(\IC_{Z^\theta})|_{X^\theta-\Delta_X}\simeq 
\underset{{\mathfrak P}(\theta),\,|{\mathfrak P}(\theta)|\neq 1}\oplus\, 
i_{{\mathfrak P}(\theta)}{}_*(\IC^{{\mathfrak P}(\theta)}\otimes \pi_M^*(\E_{{\mathfrak P}(\theta)}))
[-|{\mathfrak P}(\theta)|]|_{X^\theta-\Delta_X},$$
where $\E_{{\mathfrak P}(\theta)}$ is an order $2$ local system on the preimage of $X^\theta-\Delta_X$ in 
$X^{{\mathfrak P}(\theta)}$,
which can be non-trivial only for ${\mathfrak P}(\theta)$ of the form specified above. 
Of course, later we will have to show that all 
the $\E_{{\mathfrak P}(\theta)}$'s are necessarily trivial.

By combining this with \corref{semisimple} we obtain:

\begin{equation}    \label{directsum}
\s^{\theta}{}^!(\IC_{Z^\theta})\simeq \underset{{\mathfrak P}(\theta),\,|{\mathfrak P}(\theta)|\neq 1}\oplus\, 
i_{{\mathfrak P}(\theta)}{}_*(\K^{{\mathfrak P}(\theta)})[-|{\mathfrak P}(\theta)|]\oplus \K^\theta,
\end{equation}
where $\K^\theta$ is a complex supported on $\Mod_M^{+,\theta}|_{\Delta_X}\simeq \Gr_{M,X}^{+,\theta}$ and
$\K^{{\mathfrak P}(\theta)}$ is the Goresky-MacPherson extension of 
$\IC^{{\mathfrak P}(\theta)}|_{X^\theta-\Delta_X}\otimes \pi_M^*(\E_{{\mathfrak P}(\theta)})$ to the whole
$\Gr_M^{{\mathfrak P}(\theta)}$.

\medskip

To prove the theorem we have to understand the complex $\K^\theta$. This will be done by analyzing
the *-restriction of the LHS of \eqref{directsum} to $\Mod_M^{+,\theta}|_{\Delta_X}$.

\ssec{}

Recall that for $x\in X$ we denoted by $\Delta_x$ the embedding 
$\on{pt}\hookrightarrow X\overset{\Delta_X}\hookrightarrow X^\theta$. To simplify the notation, instead of
$\s^{\theta}{}^!(\IC_{Z^\theta})|_{\Delta_X}$ we will compute $\s^{\theta}{}^!(\IC_{Z^\theta})|_{\Delta_x}$.
We will prove the following assertion:

\begin{prop}  \label{upperestimate}
$\s^{\theta}{}^!(\IC_{Z^\theta})|_{\Delta_x}$, being a complex of sheaves on $\Mod_M^{+,\theta}|_{\Delta_x}\simeq
\Gr_M^{+,\theta}$,
is concentrated in perverse cohomological degrees $\leq 0$. Its $0$-th perverse cohomology 
can be identified with $\on{Loc}(U(\check{\mathfrak u}(P))_{\theta})$.
\end{prop}

\noindent{\it Remark.}
Note that the {\it a posteriori} proven \corref{restr2} implies that the above complex 
has perverse cohomology only in dimension $0$. 

\medskip

Recall our definition of the central fiber $\SS^\theta$. 
Using \propref{directandinverse} and base change, we obtain that
$$\s^{\theta}{}^!(\IC_{Z^\theta})|_{\Delta_x}\simeq \pi_P{}_!(\IC_{Z^\theta}|_{\SS^\theta}).$$

The following is a refinement of \propref{upperestimate}:

\begin{prop}   \label{upperestimaterefined}
Let $\theta'$ be as above.
\begin{enumerate}

\item
$\pi_P{}_!(\IC_{Z^\theta}|_{_{\theta'}\SS^\theta})$ lives in strictly negative cohomological dimensions if $\theta'\neq 0$.

\item
The complex $\pi_P{}_!(\IC_{Z^\theta}|_{_0\SS^\theta})$ lives in cohomological dimensions $\leq 0$.

\item
$h^0(\pi_P{}_!(\IC_{Z^\theta}|_{_0\SS^\theta}))\simeq \on{Loc}(U(\check{\mathfrak u}(P))_{\theta})$.

\end{enumerate}
\end{prop}

The proof of \propref{upperestimaterefined} will use the following facts about the geometry of the affine Grassmannian,
whose proofs will be given in \secref{affgrass}:

Let us denote by $\tf^\theta$ the natural map $\Gr_P^\theta\to \Gr^\theta_M$.

\begin{thm}  \label{thmaffgrass}
We have:
\begin{enumerate}

\item
$\tf^\theta_!(\Ql{}_{\Gr_P^\theta\cap \Gr_{U(P^-)}})$ as a complex of sheaves on $\Gr_M^{+,\theta}$
lives in the perverse cohomological dimensions 
$\leq \langle \theta, 2(\rhoch_G-\rhoch_M)\rangle$.

\item 
$h^{\langle \theta, 2(\rhoch_G-\rhoch_M)\rangle}(\tf^\theta_!(\Ql{}_{\Gr_P^\theta\cap \Gr_{U(P^-)}}))\simeq
\on{Loc}(U(\check{\mathfrak u}(P))_{\theta})$.
\end{enumerate}

\end{thm}

\ssec{Proof of \propref{upperestimaterefined}}

First, from \propref{identification} we can compute the dimension of $Z^\theta$ and we obtain 
$\langle \theta, 2(\rhoch_G-\rhoch_M)\rangle$. Since $_0Z^\theta$ is contained in $Z_{max}^\theta$
(and $Z_{max}^\theta$ is smooth, according to \corref{opensmooth}),
$$\IC_{Z^\theta}|_{_0\SS^\theta}\simeq \Ql_{_0\SS^\theta}[\langle \theta, 2(\rhoch_G-\rhoch_M)\rangle].$$

Therefore, points 2 and 3 of the proposition follow immediately
from \thmref{thmaffgrass} combined with \propref{central fiber}.

\medskip

To prove point 1 of the proposition, let us first take $\theta'\neq \theta$. However,
$_{\theta'}\SS^\theta$ is contained in $_{\theta'}Z^\theta$, we will be able to use the induction hypothesis
to calculate $\IC_{Z^\theta}|_{_{\theta'}\SS^\theta}$:

\medskip

Recall the identification
$$_{\theta'}\SS^\theta\simeq (\Gr_P^{\theta-\theta'}\cap \Gr_{U(P^-)})\underset{\Gr_M}\times \Conv_M^{+,\theta'}$$
of \secref{various strata}. 

Let us recall also the following construction:

Projection $pr'$ realizes the convolution diagram $\Conv_M$ as a fibration over $\Gr_M$ with the typical
fiber isomorphic to $\Gr_M$ itself. Hence, starting with a spherical perverse sheaf $\S$ on $\Gr_M$ and
an arbitrary complex $\S'$ on $\Gr_M$, we can define their twisted external product 
$\S\widetilde{\boxtimes}\S'\in\on{D}(\Conv_M)$, which is ``$\S'$ along the base'', and
``$\S$ along the fiber''. The convolution of $\S$ and $\S'$ is by definition the complex on $\Gr_M$ equal to
$pr_!(\S\widetilde{\boxtimes}\S')$. It is a basic fact (cf. \cite{Ga}) that if $\S'$ is a perverse sheaf,
then its convolution with any $\S$ as above is perverse as well.

Similarly, if $\S$ is a spherical perverse sheaf on $\Gr_M$ and $\S''$ is an arbitrary complex on 
$\Gr_P^{\theta-\theta'}\cap \Gr_{U(P^-)}$ we can construct the complex $\S\widetilde{\boxtimes}\S''$
on $$\Conv_M^{+,\theta'}\underset{\Gr_M}\times (\Gr_P^{\theta-\theta'}\cap \Gr_{U(P^-)}).$$

\medskip

By combining \thmref{restr} with \lemref{restrM} and \propref{identification} we obtain that
$\IC_{Z^\theta}|_{_{\theta'}\SS^\theta}$ is the direct sum
$$\underset{{\mathfrak P}(\theta')}\oplus\, 
(\on{Loc}(\underset{k}\otimes\, \on{Sym}^{n_k}(\check{\mathfrak u}(P)_{\theta'_k}))[2\cdot|{\mathfrak P}(\theta')|])
\widetilde{\boxtimes} \Ql[\langle\theta-\theta',2(\rhoch-\rhoch_M)\rangle].$$

Now, projection $\pi_P:{}_{\theta'}\SS^\theta\to \Gr_M^{+,\theta}$ in the above description 
of $_{\theta'}\SS^\theta$ corresponds
to $\Conv_M^{+,\theta'}\underset{\Gr_M}\times (\Gr_P^{\theta-\theta'}\cap \Gr_{U(P^-)})\to \Conv_M^{+,\theta'}
\overset{pr}\to \Gr_M$. 

Therefore, $\pi_P{}_!(\IC_{Z^\theta}|_{_{\theta'}\SS^\theta})$ is the sum over ${\mathfrak P}(\theta')$
of convolutions of 
$$\tf^{\theta-\theta'}_!(\Ql{}_{\Gr_P^{\theta-\theta'}\cap 
\Gr_{U(P^-)}})[2\cdot |{\mathfrak P}(\theta')|+\langle\theta-\theta',2(\rhoch-\rhoch_M)\rangle]$$ 
with the spherical perverse sheaf $\on{Loc}(\underset{k}\otimes\, \on{Sym}^{n_k}(\check{\mathfrak u}(P)_{\theta_k}))$. 
The important thing for us is that
$|{\mathfrak P}(\theta')|>0$: using \thmref{thmaffgrass}(2), 
we obtain that $\pi_P{}_!(\IC_{Z^\theta}|_{_{\theta'}\SS^\theta})$
lies in strictly negative cohomological degrees.

Since the convolution of a complex lying in negative perverse cohomological dimensions on $\Gr_M$
with a spherical perverse sheaf is again a complex lying in negative cohomological dimensions, point 1 of
the proposition follows for $\theta'\neq\theta$.

\medskip

Finally, let us consider $\theta'=\theta$. In this case, $\pi_P:{}_\theta\SS^\theta\to \Gr_M^{+,\theta}$
is an isomorphism and it suffices to observe, that by the very definition of intersection cohomology,
$\IC_{_\theta\SS^\theta}$ lives in strictly negative cohomological degrees.

Thus, \propref{upperestimaterefined} is proved modulo \thmref{thmaffgrass}, which will be dealt with later.

\ssec{}

Let us go back to the isomorphism of Equation \eqref{directsum}.
At this point we are ready to prove that the local systems $\E_{{\mathfrak P}(\theta)}$ are all trivial.
For that purpose, we can assume that ${\mathfrak P}(\theta)$ corresponds to
$\theta=2\cdot \theta'$, as above. Consider 
$$\grb_{M,X}^{\flat(\theta')}\times \grb_{M,X}^{\flat(\theta')} \to\Gr_M^{{\mathfrak P}(\theta)} \overset
{i_{{\mathfrak P}(\theta)}}\longrightarrow \Mod_M^{+,\theta}.$$

By induction hypothesis and \propref{factor} we have that over $X\times X-\Delta_X$,
\begin{equation}   \label{Z2monodromy}
(\s^{\theta}{}^!(\IC_{Z^\theta}))|_{\grb_{M,X}^{\flat(\theta')}\times \grb_{M,X}^{\flat(\theta')}}\simeq
\on{Loc}_X(\on{Sym}(\check{\mathfrak u}(P))_{\theta'})\boxtimes \on{Loc}_X(\on{Sym}(\check{\mathfrak u}(P))_{\theta'})[-2].
\end{equation}

The group $\ZZ_2$ acts in a natural way on $\grb_{M,X}^{\flat(\theta')}\times \grb_{M,X}^{\flat(\theta')}$
and we have to show that the $\ZZ_2$-equivariant structure on the LHS of \eqref{Z2monodromy} corresponds to
the tautological $\ZZ_2$-equivariant structure on the RHS. 

Let us apply a relative version of \propref{upperestimaterefined} for $\theta'$, in which instead of a fixed
$x\in X$ we have a pair of distinct points on $X$. We obtain an isomorphism of complexes over $X\times X-\Delta_X$

\begin{align*}
&(\s^{\theta}{}^!(\IC_{Z^\theta}))|_{\grb_{M,X}^{\flat(\theta')}\times \grb_{M,X}^{\flat(\theta')}}\simeq 
(\pi_P{}_!(\IC_{Z^\theta}))|_{\grb_{M,X}^{\flat(\theta')}\times \grb_{M,X}^{\flat(\theta')}}\simeq \\
&h^{-top}((\tf^{\theta'}\boxtimes\tf^{\theta'})_!(\Ql_{\Gr_{P,X}^{\theta'}\cap \Gr_{U(P^-),X}\times 
\Gr_{P,X}^{\theta'}\cap \Gr_{U(P^-),X}}))\simeq \\
&\on{Loc}_X(U(\check{\mathfrak u}(P))_{\theta'})\boxtimes \on{Loc}_X(U(\check{\mathfrak u}(P))_{\theta'})[-2], 
\end{align*}
where $top=2\cdot (1+\langle \theta', 2(\rhoch_G-\rhoch_M)\rangle)$ and $\Gr_{P,X}^\theta$ and $\Gr_{U(P^-),X}$ are 
the corresponding relative (over $X$) versions of  $\Gr_{P}^\theta$ and $\Gr_{U(P^-)}$, respectively.

\medskip

The last isomorphism, by construction, intertwines the natural $\ZZ_2$-structure on
$(\s^{\theta}{}^!(\IC_{Z^\theta}))|_{\grb_{M,X}^{\flat(\theta')}\times \grb_{M,X}^{\flat(\theta')}}$ and the tautological
$\ZZ_2$-structure on the external product
$\on{Loc}(U(\check{\mathfrak u}(P))_{\theta'})\boxtimes \on{Loc}(U(\check{\mathfrak u}(P))_{\theta'})$.
By comparing with \eqref{Z2monodromy} we obtain the required assertion.

\ssec{}

To prove the theorem over $\Fqb$ it remains to analyze the term $\K^\theta$. Note that there is not more than one
${\mathfrak P}(\theta)$ with $|{\mathfrak P}(\theta)|=1$. We will denote it by ${\mathfrak P}^0(\theta)$.
By definition, $\Gr_M^{{\mathfrak P}^0(\theta)}\simeq \grb_{M,X}^{\flat(\theta)}$.

We have to show that
$$\K^\theta\simeq i_{{\mathfrak P}^0(\theta)}{}_*(\on{Loc}_X(\check{\mathfrak u}(P)_{\theta})) \otimes 
(\Ql(\frac{1}{2})[1])^{-1}.$$

By the definition of IC, since $\K^\theta|_{\Delta_X}$ is a direct summand of $\s^{\theta}{}^!(\IC_{Z^\theta})$, it can
have perverse cohomology only in degrees $\geq 1$. Let us now restrict both sides of \eqref{directsum} to 
$\Gr^{+,\theta}_{M,X}\simeq \Mod_M^{+,\theta}|_{\Delta_X}$ and apply the cohomological truncation $\tau^{\geq 1}$.
Using \lemref{restrM} on the one hand, and the relative (over $X$) version of \propref{upperestimate} on the other hand,
we obtain:
$$\on{Loc}_X((U(\check{\mathfrak u}(P)))_\theta)[-1]\simeq \underset{{\mathfrak P}(\theta),\,|{\mathfrak P}(\theta)|\neq 1}\oplus\, 
\on{Loc}_X(\underset{k}\otimes\, \on{Sym}^{n_k}(\check{\mathfrak u}(P)_{\theta_k}))[-1]\oplus \K^\theta|_{\Delta_X}.$$

Hence, $\K^\theta[1]$ is a perverse sheaf. Moreover, since $U(\check{\mathfrak u}(P))$ and $\on{Sym}(\check{\mathfrak u}(P))$
are (non-canonically) isomorphic as $\check M$-modules, the comparison of multiplicities forces 
$\K^\theta[1]|_{\Delta_X}\simeq \on{Loc}_X(\check{\mathfrak u}(P)_\theta)$.

\ssec{}  \label{Weil}

Now let us restore the $\Fq$-structure on $\s^{\theta}{}^!(\IC_{Z^\theta})$. To complete the proof of the theorem,
by induction, it suffices to show that the arrow
$$\K^\theta|_{\Delta_X}\simeq \on{Loc}_X(\check{\mathfrak u}(P)_\theta)\otimes (\Ql[1](\frac{1}{2}))^{\otimes -1}\to
\s^{\theta}{}^!(\IC_{Z^\theta}),$$ which is known to split over $\Fqb$, splits over $\Fq$ as well. For that, it is enough to show
that the complex $\s^{\theta}{}^!(\IC_{Z^\theta})|_{\Delta_X}$ is semisimple.

We know already that $\s^{\theta}{}^!(\IC_{Z^\theta})|_{\Delta_X}\otimes (\Ql[1](\frac{1}{2}))$ 
has perverse cohomology only in
dimension $0$, which is equal, as in \propref{upperestimaterefined}, to 
$$h^{0}(\tf^\theta_!(\Ql{}_{\Gr_{P,X}^{\theta}\cap \Gr_{U(P^-),X}})\otimes 
(\Ql[1](\frac{1}{2}))^{\langle \theta, 2(\rhoch_G-\rhoch_M)\rangle}).$$

The needed result follows from the fact that the isomorphism of \thmref{thmaffgrass}(3) 
is compatible with the $\Fq$-structure in the sense that
$$h^{0}(\tf^\theta_!(\Ql{}_{\Gr_{P,X}^{\theta}\cap \Gr_{U(P^-),X}})\otimes 
(\Ql[1](\frac{1}{2}))^{\langle \theta, 2(\rhoch_G-\rhoch_M)\rangle})\simeq 
\on{Loc}_X((U(\check{\mathfrak u}(P)))_\theta).$$

\section{Intersections of semi-infinite orbits in the affine Grassmannian}  \label{affgrass}

\ssec{The restriction functors}

Let $\O_x$ (resp., $\K_x$) denote the completed local ring (resp., local field) at $x$.
We can form the group-schemes $G(\O_x)$, $P(\O_x)$, $U(P)(\O_x)$ and the corresponding
group-ind-schemes $G(\K_x)$, $P(\K_x)$, $U(P)(\K_x)$. Note, however, that the latter
is not only a group-ind-scheme, but also an ind-group-scheme, i.e. an inductive limit
of group-schemes.

Let $\nu\in\Lambda$ be $M$-dominant and let $\theta$ be its image under 
$\Lambda\to \Lambda_{G,P}$. Let us denote by
$\Gr_P^{\nu}$ the pre-image $(\tf_P^{\theta})^{-1}(\Gr^\nu_M)
\subset \Gr_P^{\theta}$. The schemes $\Gr_P^{\nu}$ are nothing but orbits of the group
$U(P)(\K_x)\cdot M(\O_x)$ on $\Gr_G$. We will denote by $\tf^\nu_P$ the
restriction of $\tf_P^{\theta}$ to $\Gr_P^{\nu}$.

The goal of this section is to prove \thmref{thmaffgrass}. The starting point
is the following result, which describes the intersections of $\Gr_P^\theta$ with
$\Gr^\lambda_G$ inside the affine Grassmannian $\Gr_G$ (cf. \cite{BD}, \cite{BG} and \cite{MV}).

For a $G$-dominant (resp., $M$-dominant) coweight $\lambda$, let $V^\lambda$ (resp., $V^\lambda_M$)
denote the corresponding irreducible representation of $\check G$ (resp., $\check M$).

\begin{thm}      \label{restriction}
Let $\lambda$ be a dominant integral coweight of $G$.
\begin{enumerate}
\item The intersection $\Gr_P^\nu\cap\Gr_G^\lambda$ has dimension
$\leq \langle\nu+\lambda,\rhoch_G\rangle$. 
\item The irreducible components of $\Gr_P^\nu\cap\Gr_G^\lambda$ of dimension $\langle\nu+\lambda,\rhoch_G\rangle$
form a basis for $\Hom_{\check M}(V^\nu_M,\on{Res}^{\check G}_{\check M}(V^\lambda))$.
\end{enumerate}
\end{thm}

\ssec{The case $P=B$} 

We will first consider the situation when $P=B$. Note that in our notation
$\Bun_{B^-}^0$ is the same as $\Bun_{U(B^-)}$.

In this case,
$\Lambda_{G,P}=\Lambda$ and for two elements $\nu,\mu\in\Lambda$ let us
consider the intersection $\Gr^{\nu-\mu}_{B}\cap \Gr^{-\mu}_{B^-}$. 

\medskip

First, it is easy to see that the action of $t^{\mu}\in T(\K_x)$ on $\Gr_G$ identifies 
$\Gr^\nu_{B}\cap \Gr^0_{B^-}$ with $\Gr^{\nu-\mu}_{B}\cap \Gr^{-\mu}_{B^-}$
for any $\mu\in\Lambda$. To prove the theorem, it suffices to show that for a given 
$\nu\in \Lambda^{pos}$ and some $\mu\in \Lambda$, the intersection
$\Gr^{\nu-\mu}_{B}\cap \Gr^{-\mu}_{B^-}$ is of dimension $\leq \langle \nu,\rhoch_G\rangle$
and
$$H_c^{2\langle\nu,\rhoch\rangle}(\Gr_B^{\nu-\mu}\cap\Gr_{B^-}^{-\mu})\simeq
U(\check{\mathfrak u})_{\nu}.$$

\begin{prop}   \label{contained}
For a fixed $\nu$ and $\mu$ deep enough in the dominant
chamber, the intersection $\Gr^{\nu-\mu}_{B}\cap \Gr^{-\mu}_{B^-}$
is contained inside $\Gr_G^{-w_0(\mu)}$.
\end{prop}

\begin{proof}

Let us identify $\Gr^{0}_{B^-}$ with the quotient
$U(B^-)(\K_x)/U(B^-)(\O_x)$.

Since we know already that $\Gr^{\nu}_{B}\cap \Gr^{0}_{B^-}$ is a scheme
of finite type, for $\mu$ deep enough in the dominant chamber,
the preimage of $\Gr^{\nu}_{B}\cap \Gr^{0}_{B^-}$ under the projection
$U(B^-)(\K_x)\to U(B^-)(\K_x)/U(B^-)(\O_x)$ is contained inside the subgroup
$$\on{Ad}_{t^{\mu}}(U(B^-)(\O_x))\subset U(B^-)(\K_x).$$

\medskip

Let us now consider $\Gr^{\nu-\mu}_{B}\cap \Gr^{-\mu}_{B^-}$, which via
the action of $t^{\mu}$ can be identified with $\Gr^{\nu}_{B}\cap \Gr^{0}_{B^-}$.

We can view $\Gr^{-\mu}_{B^-}$ as a quotient $U(B^-)(\K_x)/\on{Ad}_{t^{-\mu}}(U(B^-)(\O_x))$,
via the action of $U(B^-)(\K_x)$ on $t^{\mu}$, viewed as an element of $\Gr_T\subset \Gr^{-\mu}_{B^-}\subset \Gr_G$.
We obtain that the preimage of $\Gr^{\nu-\mu}_{B}\cap \Gr^{-\mu}_{B^-}$
in $U(B^-)(\K_x)$ is contained in $U(B^-)(\O_x)$.

Hence, 
$$\Gr^{\nu-\mu}_{B}\cap \Gr^{-\mu}_{B^-}\subset U(B^-)(\O_x)\cdot t^{-\mu}\subset 
G(\O_x)\cdot t^{-\mu}= \Gr_G^{-w_0(\mu)}.$$

\end{proof}

The above proposition implies the dimension estimate 
$\dim(\Gr^{\nu-\mu}_{B}\cap \Gr^{-\mu}_{B^-})\leq \langle \nu-\mu,\rhoch_G\rangle$
right away.

Indeed, we may assume that $\mu$ is such that $\Gr^{\nu-\mu}_{B}\cap \Gr^{-\mu}_{B^-}\subset \Gr_G^{-w_0(\mu)}$.
However, \thmref{restriction}(1) implies that
$\dim(\Gr_G^{-w_0(\mu)}\cap \Gr^{\nu-\mu}_{B})\leq \langle \nu,\rhoch_G\rangle$.

\medskip

To prove the other statements of the theorem, observe that for $\mu$ as above the irreducible components
of $\Gr^{\nu-\mu}_{B}\cap \Gr^{-\mu}_{B^-}$ of dimension $\langle \nu,\rhoch_G\rangle$ are naturally
a subset among the irreducible components of $\Gr_G^{-w_0(\mu)}\cap \Gr^{\nu-\mu}_{B}$ of
the same dimension. 

Let us show that the generic point of
every irreducible component $K$ of the intersection $\Gr_G^{-w_0(\mu)}\cap \Gr^{\nu-\mu}_{B}$
of dimension $\langle \nu,\rhoch_G\rangle$ is contained in $\Gr^{\nu-\mu}_{B}\cap \Gr^{-\mu}_{B^-}$.

\medskip

Suppose the contrary. Then there exists $\mu'\in\Lambda$, such that 
the generic point of $K$ is contained in $\Gr_G^{-w_0(\mu)}\cap \Gr^{-\mu'}_{B^-}$. 
However, it is easy to see that $\Gr_G^{-w_0(\mu)}\cap \Gr^{-\mu'}_{B^-}\neq \emptyset$ implies
$\mu-\mu'\in\Lambda^{\on{pos}}$. 

However, as we have shown above,
$\dim(\Gr^{\nu-\mu}_{B} \cap \Gr^{-\mu'}_{B^-})\leq \langle \nu+\mu'-\mu,\rhoch_G\rangle$,
which is smaller than the dimension of $K$.

\medskip

Thus, we obtain that
$$H_c^{2\langle\nu,\rhoch\rangle}(\Gr_B^{\nu-\mu}\cap\Gr_{B^-}^{-\mu})\simeq
H_c^{2\langle\nu,\rhoch\rangle}(\Gr_G^{-w_0(\mu)}\cap \Gr^{\nu-\mu}_{B}).$$
However, according to \thmref{restriction}(2), the RHS of the above equation 
can be canonically identified with the $\nu-\mu$-weight
space in the irreducible $\check G$-representation with highest weight $-w_0(\mu)$.

The latter, when $\mu$ is large compared to $\nu$, 
is isomorphic to $U(\check{\mathfrak u})_{\nu}$
via the action on the lowest weight vector.

\ssec{The general case}

We fix $\theta$ and $\nu\in \Lambda$ such that $\Gr_M^{\nu}\subset \Gr_M^{+,\theta}$.
Since $\Gr_M^\nu$ is simply-connected, it suffices to show that 
each intersection $\Gr_P^\nu\cap \Gr_{U(P^-)}$ is of dimension
$\leq \langle \nu, \rhoch_G\rangle$ 
and that the number of its irreducible components of
dimension exactly $\langle \nu, \rhoch_G\rangle$
equals the dimension of $\Hom_{\check M}(V_M^\nu,U(\check{\mathfrak u}(P)))$.

For an $M$-dominant weight $\mu$ let us consider the corresponding $\Gr_{P^-}^\mu\subset \Gr_G$.
Note that for $\mu=0$ this sub-scheme coincides with $\Gr_{U(P^-)}$.

\medskip

Let $\Lambda'_{G,P}\subset \Lambda_G$ denote the lattice of cocharacters of
the center $Z(M)$ of $M$. If $\mu'\in \Lambda'_{G,P}$, the action of the corresponding $t^{\mu'}\in Z(M)(\K_x)$ identifies
$\Gr_P^{\nu}\cap \Gr_{P^-}^{\mu}$ with $\Gr_P^{\nu-\mu'}\cap \Gr_{P^-}^{\mu-\mu'}$.

\begin{prop} \label{contained non-prince}
Let $\mu'\in \Lambda'_{G,P}$ be $G$-dominant and deep enough on the corresponding
wall of the Weyl chamber. Then the intersection
$\Gr_P^{\nu-\mu'}\cap \Gr_{P^-}^{\mu-\mu'}$ is contained in $\Gr_G^{w_0(w^M_0(\mu-\mu'))}$.
\end{prop}

\begin{proof}

The initial observation is that each $\Gr_P^{\nu}\cap \Gr_{P^-}^{\mu}$
is a scheme of finite type. We know this fact for $\mu=0$, since the above
intersection is a locally closed sub-scheme in the Zastava space $Z^\theta$.

In general, this assertion can be proven either by introducing the
corresponding analog of the Zastava space over a global curve, or
by a straightforward local argument.

\medskip

Let us view $\Gr_M^{\mu}$ as a sub-scheme of $\Gr^\mu_{P^-}$, such that 
$\Gr_{P^-}^{\mu}=U(P^-)(\K_x)\cdot \Gr_M^{\mu}$. 

As in the case $P=B$, we obtain that that the preimage of 
$\Gr_P^\nu\cap \Gr_{P^-}^{\mu}$ under
$$U(P^-)(\K_x)\times \Gr_M^{\mu}\to \Gr_{P^-}^{\mu}$$
is contained in a sub-scheme of the form $\on{Ad}_{t^{\mu'}}(U(P^-)(\O_x))\times \Gr_M^{\mu-\mu'}$.

Hence, the action of $t^{-\mu'}$ maps $\Gr_P^\nu\cap \Gr_{P^-}^{\mu}$ inside
$$U(P^-)(\O_x)\cdot \Gr_M^{\mu-\mu'}\subset \Gr_G^{w_0(w_0^M(\mu-\mu'))}.$$

\end{proof}

The rest of the proof is similar to the case of $P=B$:

\medskip

From the above proposition we obtain that the intersection $\Gr_P^{\nu}\cap \Gr_{P^-}^{\mu}$
is of dimension $\leq\langle \nu-w^M_0(\mu), \rhoch_G\rangle$. In particular,
$\Gr_P^\nu\cap \Gr_{U(P^-)}$ is of dimension $\leq \langle \nu, \rhoch_G\rangle$.

Moreover, as in the previous case, we obtain that there is a bijection between the set
of irreducible components of $\Gr_P^\nu\cap \Gr_{U(P^-)}$ of dimension $\langle \nu, \rhoch_G\rangle$
and the set of irreducible components of the same dimension of $\Gr_P^{\nu-\mu'}\cap \Gr_G^{-w_0(\mu')}$,
where $\mu'$ is large enough.

However, \thmref{restriction}(2) implies that the latter set parameterizes a basis of
$$\Hom_{\check M}(V_M^{\nu-\mu'},V_G^{-w_0(\mu')}),$$ which, since $\mu'$ is large
compared to $\nu$, can be identified with $\Hom_{\check M}(V_M^\nu,U(\check{\mathfrak u}(P)))$.

\section{Intersection cohomology of $\BunPb$}

In this section we will be concerned with describing explicitly the intersection cohomology
sheaf on $\BunPb$. First, we introduce an analogue, $Z^{nv,\theta}$, of the Zastava spaces for $\BunPb$
(here the superscript $nv$ stands for ``naive'').

\medskip

By definition, $Z^{nv,\theta}$ is a scheme classifying
the data of $(x^\theta,\F_G,\beta)$, where $x^\theta\in X^\theta$, $\F_G$ is a $G$-bundle and
$\beta$ is a trivialization of $\F_G$ off the support of $x^\theta$, such that 
for every $G$-dominant weight $\lambdach$ orthogonal to $\on{Span}(\alpha_\i)$, $\i\in\I_M$
the induced meromorphic maps
$$\L_{\F^0_{M/[M,M]}}^{\lambdach}\to \V_{\F_G}^{\lambdach} \text{ and }
\V_{\F_G}^{\lambdach}\to \L_{\F^0_{M/[M,M]}}^{\lambdach}$$
induce a regular map
$\L_{\F^0_{M/[M,M](-x^\theta)}}^{\lambdach}\to \V_{\F_G}^{\lambdach}$
and a regular and {\it surjective} map 
$\V_{\F_G}^{\lambdach}\to \L_{\F^0_{M/[M,M]}}^{\lambdach}$.

\medskip

There is a natural proper map $Z^\theta\to Z^{nv,\theta}$, which corresponds to
``forgetting'' the data of $(\F_M,\beta_M)$. In addition, $Z^{nv,\theta}$ contains an
open sub-scheme $Z^{nv,\theta}_{max}$ corresponding to the locus, where the maps
$\L_{\F^0_{M/[M,M](-x^\theta)}}^{\lambdach}\to \V_{\F_G}^{\lambdach}$ are maximal embeddings,
over which we have an isomorphism $Z^\theta_{max}\to  Z^{nv,\theta}_{max}$.

\medskip

As in the case of $Z^\theta$, one easily establishes the {\it Factorization property}
for $Z^{nv,\theta}$:

\begin{equation} \label{naive factorization}
(X^{\theta_1}\times X^{\theta_2})_{disj}\underset{X^\theta}\times Z^{nv,\theta}\simeq 
(X^{\theta_1}\times X^{\theta_2})_{disj}\underset{X^{\theta_1}\times X^{\theta_2}} 
\times (Z^{nv,\theta_1}\times Z^{nv,\theta_2}).
\end{equation}

\medskip

Finally, the spaces $Z^{nv,\theta}$ model the 
singularities of $\BunPb$ in the same
sense as $Z^\theta$ models the singularities of $\BunPbw$. In other words, an analogue
of \propref{identification} holds, whose formulation we leave to the reader.

\ssec{}

It turns out, that although the stack $\BunPb$ is ``simpler'' than $\BunPbw$,
the description of its intersection cohomology sheaf is
more involved, and in particular it relies on the description of $\IC_{\BunPbw}$.

To formulate the main theorem, we introduce the following notation:

\medskip

Recall that the functors $H^\bullet$ and $\on{Loc}$ establish the quasi-inverse
equivalences between the category of spherical perverse sheaves on $\Gr_M$ and
the category of $\check M$-modules.

Under this equivalence, the multiplication
by the first Chern class $c_1(\det)$ of the determinant line bundle on 
$\Gr_M$ corresponds to the action of a principal nilpotent element 
$e\in\Lie(\check M)$. Moreover, the cohomological grading of $H^\bullet$
corresponds to the action of a semisimple $h\in\Lie(\check M)$ contained in
a uniquely defined principal $\mathfrak{sl}_2$-triple $(e,f,h)$ in
$\Lie(\check M)$. For a $\check M$-module $V$ the 
$\ZZ$-grading arising from the action of $h$ 
is given by the following rule: the weight subspace $V_\eta$
has degree $\langle\eta,2\rhoch_M\rangle$. This $\ZZ$-grading on $V$ will
be called {\em the principal grading}. 

For $V$ as above, we will denote by $V^f$ the subspace annihilated by $f$. We will
consider it as a graded vector space, via the principal grading.

\medskip

We define the functor $\ol{\on{Loc}}\,(V)$ from the category of $\ZZ$-graded vector spaces to complexes
over $\on{Spec}(\Fq)$ by setting 
$$\ol{\on{Loc}}\,(V)=\underset{n}\oplus\, V_n\otimes [-n](-\frac{n}{2}).$$
In particular, we will apply the functor $\ol{\on{Loc}}$ to $\check M$-modules $V$ (or their
direct summands, such as $V^f$), endowed with the principal grading.

\begin{thm}  \label{IC of BunPb on theta}
The restriction of $\IC_\BunPb$ to $$_\theta\BunPb\simeq X^\theta\times \Bun_P\simeq \underset{k}\Pi\, X^{(n_k)}\times \Bun_P$$
can be identified with the direct sum over $\fP(\theta)$
of the direct images under 
$X^{\fP(\theta)}\times\Bun_P\to X^\theta\times\Bun_P$ of
$$\IC_{\Bun_P}\boxtimes \left(\underset{k}\boxtimes\, \left(\ol{\on{Loc}}\,(\frn^f_{\theta_k})\right)^{(n_k)} \right) \otimes
\left(\Ql\left(\frac{1}{2}\right)[1]\right)^{\otimes 2\cdot |\fP(\theta)|},$$
where each $\ol{\on{Loc}}\,(\frn_{\theta_k})$ is viewed as a constant local system on $X$ and the superscript
$(n_k)$ designates the $n_k$-th symmetric power.
\end{thm}

As a corollary, we obtain the following description of the restriction of $\IC_{\BunPb}$ to the strata
$_{\fA(\theta)}\BunPb$:

\begin{thm} \label{IC of BunPb on A(theta)}
Let $\fA(\theta)$ be a partition $\theta=\underset{k}\sum\, n_k\cdot\theta_k$. 
The $*$-restriction of $\IC_\BunPb$ to  
$_{\fA(\theta)}\BunPb\simeq \oX^{\fA(\theta)}\times \Bun_P$ 
is isomorphic to $\IC_{_{\fA(\theta)}\BunPb}$ tensored by the 
complex\footnote{In the published version~\cite{BFGM} there is a misprint
in the formula for the $*$-restriction: instead of the symmetric $n_k$-th
power it contains the tensor $n_k$-th power, that is 
$$\bigotimes_k\left(\underset{i\geq 0}\oplus \,\ol{\on{Loc}}\left(\Sym^i(\frn^f)_{\theta_k}\right)
(i)[2i]\right)^{\otimes n_k}\otimes\left(\Ql\left(\frac{1}{2}\right)[1]\right)^{\otimes-|\fA(\theta)|}.$$}
$$\bigotimes_k\left(\underset{i\geq 0}\oplus \,\ol{\on{Loc}}\left(\Sym^i(\frn^f)_{\theta_k}\right)
(i)[2i]\right)^{(n_k)}\otimes\left(\Ql\left(\frac{1}{2}\right)[1]\right)^{\otimes-|\fA(\theta)|}.$$
\end{thm}

\noindent{\em Remark.} Suppose $G=SL_n$. Then the parabolic subgroups of $G$
are numbered by the ordered partitions $n=n_1+\ldots+n_k,\ n_i>0$. Suppose
a parabolic subgroup $P$ corresponds to a non-decreasing partition
$n=n_1+\ldots+n_k,\ 0<n_1\leq n_2\leq\ldots\leq n_k$. In this case 
\thmref{IC of BunPb on A(theta)} was proved by A.~Kuznetsov in the Summer 1997 (unpublished).
His proof made use of {\em Laumon's compactification} $\Bun_P^L$ ~\cite{La}.
Namely, $\Bun_P^L$ is always smooth (see {\em loc. cit.}) and equipped with
a natural dominant representable projective morphism $\varpi:\ \Bun_P^L\to
\BunPb$. In case $P$ corresponds to a non-decreasing partition, A.~Kuznetsov
proved that $\varpi$ is small, and computed the cohomology of its fibers.

Let us mention that in case $G=SP(4)$, and $P$ corresponding to the Dynkin
sub-diagram formed by the {\em long} simple root, $\BunPb$ {\em does not
admit} a small resolution, as can be seen from the calculation of IC stalks
in codimension 5 (the existence of such resolution would imply that a fiber
has cohomology $\Ql\oplus\Ql[-4]$).

\ssec{} 

We will deduce \thmref{IC of BunPb on theta} from \corref{restr1}. Let $\Q$ denote the direct image
of $\IC_\BunPbw$ under $\r:\BunPbw\to\BunPb$. On the one hand, from \corref{restr1} 
and Lusztig's computation ~\cite{Lu} of global cohomology of perverse sheaves on affine Grassmannians,
we obtain:

\begin{cor}  \label{restr on theta}
The $*$-restriction of $\Q$ to $_\theta\BunPb\simeq X^\theta\times \Bun_P$
is isomorphic to the direct sum over $\fP(\theta)$
of the direct images under $X^{\fP(\theta)}\times\Bun_P\to X^\theta\times\Bun_P$ of
$$\IC_{\Bun_P}\boxtimes \left(\underset{k}\boxtimes\, \left(\ol{\on{Loc}}\,(\frn_{\theta_k})\right)^{(n_k)} \right) \otimes
\left(\Ql\left(\frac{1}{2}\right)[1]\right)^{\otimes 2\cdot |\fP(\theta)|}.$$
\end{cor}

On the other hand, by the Decomposition theorem, $\Q$ is a pure complex, which contains $\IC_{\BunPb}$ as a direct
summand. Therefore, \thmref{IC of BunPb on theta} amounts to identifying the corresponding direct summand
in the formula for $\Q|_{_\theta\BunPb}$ of the above corollary. In particular, 
we obtain that $\IC_{\BunPb}|_{_\theta\BunPb}$ is a pure complex.

\medskip

Consider the main diagonal $X\times \Bun_P\to X^\theta\times \Bun_P$, which corresponds to the
partition $\fA(\theta)=\fA^0(\theta)$. Let $\S_\theta$ be the direct summand of $\IC_{\BunPb}|_{_\theta\BunPb}$,
supported on $X\times \Bun_P$. By induction and the Factorization property, 
it suffices to show that
$$\S_\theta\simeq \IC_{X\times \Bun_P}\otimes\ol{\on{Loc}}\,(\frn^f_\theta)\otimes\Ql(\frac{1}{2})[1].$$

\ssec{}

We proceed as follows:

Let $\overline{_{\fP(\theta)}\BunPb}$ denote the closure of the image of 
$X^{\fP(\theta)}\times \Bun_P\to X^\theta\times \Bun_P\simeq {}_\theta\BunPb$ in $\BunPb$.
By \corref{restr on theta} and the Decomposition theorem, we have:
$$\Q=\underset{\theta\in\Lambda^\pos_{G,P},\, \fP(\theta)}\bigoplus \Q_{\fP(\theta)},$$ 
where each $\Q_{\fP(\theta)}$ is a complex on $\overline{_{\fA(\theta)}\BunPb}$.

In particular, we obtain that
$$\Q|_{_\theta\BunPb}\simeq \underset{\theta\in\Lambda^\pos_{G,P},\, \fP(\theta)}\bigoplus
\Q_{\fP(\theta)}|_{_\theta\BunPb}.$$

\begin{lem}
For 
$0<\theta'<\theta$, none of the $\Q_{\fP(\theta')}|_{_\theta\BunPb}$ has a direct summand supported
on the main diagonal $X\times \Bun_P\subset {}_\theta\BunPb$.
\end{lem}

\begin{proof}

First, from \corref{restr1}, it is easy to see that each $\Q_{\fP(\theta')}$ has the following form:
it is the intersection cohomology sheaf of $\overline{_{\fP(\theta')}\BunPb}$
tensored with a complex over $\on{Spec}(\Fq)$.

\medskip

There is a finite map $X^{\fP(\theta')}\times \BunPb\to \BunPb$, defined as in \propref{naivestrata},
which normalizes $\overline{_{\fP(\theta')}\BunPb}$. Hence, it suffices to analyze the *-restriction to
$_\theta\BunPb$ of the direct image of $\IC_{X^{\fP(\theta')}\times \BunPb}$ under this map.

However, the preimage of $_\theta\BunPb$ in $X^{\fP(\theta')}\times \BunPb\to \BunPb$ is
$X^{\fP(\theta')}\times {}_{\theta-\theta'}\BunPb$ and we can assume that the *-restriction 
of $\IC_{X^{\fP(\theta')}\times \BunPb}$ to this sub-stack is known by induction. 

In particular, all its direct summands
are supported on sub-stacks of the form $$X^{\fP(\theta')}\times X^{\fP(\theta-\theta')}\times \Bun_P.$$
Since, $\theta'\neq 0$ and $\theta-\theta'\neq 0$, none of these sub-stacks maps onto
the main diagonal in $_\theta\BunPb$.

\end{proof}

Thus, from \corref{restr on theta}, we obtain that
$$\IC_{X\times \Bun_P}\otimes\left(\ol{\on{Loc}}\,(\frn_\theta)\right)\otimes \Ql(\frac{1}{2})[1]\simeq
\S_\theta\oplus \Q_{\fP^0(\theta)}.$$

Hence, it suffices to see that if we decompose $\frn_\theta$ as
$$\frn_\theta=\frn^f_\theta\oplus \on{Im}(e:\frn\to \frn)_\theta,$$
then the induced map
$$\S_\theta\to 
\IC_{X\times \Bun_P}\otimes\left(\ol{\on{Loc}}\,(\frn_\theta)\right)\otimes \Ql(\frac{1}{2})[1]\to
\IC_{X\times \Bun_P}\otimes\left(\ol{\on{Loc}}\,(\frn^f_\theta)\right)\otimes \Ql(\frac{1}{2})[1]$$
is an isomorphism.

\smallskip

The latter is established as follows:

\bigskip

Consider the line bundle $L$ on $\BunPbw$ equal to 
the ratio of the pull-backs of the determinant
line bundles under the maps $\BunPbw\to \Bun_G$ 
and $\BunPbw\to \Bun_M$, respectively. Its restriction to the fibers of 
$\mathfrak{r}_P:\ \BunPbw\to\BunPb$ over $_{\fA^0(\theta)}\BunPb$ is equal
to a positive power of the determinant line bundle $\det$ on $\Gr_M$. 
Hence $L$ is a relatively ample line bundle for $\BunPbw$ over $\BunPb$.
The relative hard Lefschetz theorem ~\cite{BBD} asserts that the multiplication
by $c_1(L)^i$ induces an isomorphism from $\Q^{-i}$ to $\Q^i(i)$ where
$\Q^i$ denotes the direct summand of $\Q$ in perverse cohomological degree $i$.

Let us restrict the action of $c_1(L)$ to the direct summand of 
$\Q|_{_\theta\BunPb}$ supported on the main diagonal $X\times\Bun_P$.
Under identification of this summand with
$$\IC_{X\times \Bun_P}\otimes\left(\ol{\on{Loc}}\,
(\frn_\theta)\right)\otimes \Ql(\frac{1}{2})[1]$$ this action coincides,
up to a scalar, with the action of $e$, by the very definition. 
Let us disregard Tate twists and view the above direct summand as 
a semisimple graded perverse sheaf. 

\medskip

We have the following general lemma:

\begin{lem} Let $A^\bullet$ be a graded semisimple object of an abelian
category, equipped with an endomorphism $e:\ A^\bullet\to A^{\bullet+2}$
such that $e^i:\ A^{-i}\to A^i$ is an isomorphism. Suppose $A[1]=B\oplus C$
where $B$ is concentrated in negative degrees, and $e^k:\ C^{-k}\to C^k$
is an isomorphism for any $k\geq0$. Then

a) There is a unique endomorphism $f:\ A^\bullet\to A^{\bullet-2}$ satisfying
the relations of $\mathfrak{sl}_2$ together with $e,h$ where $h|_{A^i}=i$;

b) $C=\on{Im}(e)$, and the projection $B\to A[1]/C=A[1]/\on{Im}(e)
\widetilde{\leftarrow}\on{Ker}(f)$ identifies $B$ with $\on{Ker}(f)$.
\end{lem}

The proof of \thmref{IC of BunPb on theta} is concluded by applying this lemma to
$$A^{\bullet}=\IC_{X\times \Bun_P}\otimes\left(\ol{\on{Loc}}\,(\frn_\theta)\right),\,\,
B=\S_\theta \text{ and }C=\Q_{\fP^0(\theta)}.$$

\section{Erratum}

\subsection{}
\label{(1)}
The following must be added at the end of Section 3:

The above shows that a point of $\BunPbw$ belonging 
to a stratum $_{\theta}\BunPbw$ is (locally in the
smooth topology) equivalent to a point of $_{\theta}Z^{\theta'}$
for some $\theta'$ with $\theta'-\theta\in \Lambda^{\on{pos}}_{G,P}$.
In fact, we claim that one can assume that $\theta'=\theta$.
This is needed in order to obtain the description of
$\on{IC}_{\BunPbw}$ from Theorem 4.5.

Indeed, let our point $(\F_G,\F_M,\kappa)\in {}_{\theta}\BunPbw$
have singularities at $x_1,...,x_n\in X$,
and let $\F_P$ be the corresponding $P$-bundle obtained
by saturating $\kappa$. We can find a reduction $\F_{P^-}$
of $\F_G$ to $P^-$, which is in generic position with $\F_P$ 
at $x_1,...,x_n$, and such that the induced $M$-bundle $\F'_M$
belongs to $\Bun_M^r$. Then, according to Proposition 3.2, 
$(\F_G,\F_M,\kappa,\F_{P^-})$ defines a point of 
$_{\theta}Z^{\theta'}$ for some $\theta'$. However, by the 
Factorization property, Proposition 2.4, this point is \'etale-locally 
equivalent to a point in the product $_{\theta}Z^\theta\times 
Z^{\theta'-\theta}_{\on{max}}$, and we know that the second factor
is smooth.

\subsection{}
\label{(2)}
The first paragraph of Sect. 4.1 must read as follows:

For $\theta\in \Lambda^{pos}_{G;P}$, let ${\mathfrak P}(\theta)$
denote an element of the set of partitions of $\theta$ as a sum
$\theta=\underset{k}\Sigma\, \theta_k$, where each $\theta_k$ is
such that the isotypic component $\check {\mathfrak u}(P)_{\theta_k}$
is non-trivial. Recall that in this case, 
$\check {\mathfrak u}(P)_{\theta_k}$ is irreducible as a representation
of the dual Levi $\check M$, and we will denote by 
$\sharp(\theta_k)\in \Lambda^+_M\subset \Lambda_G$ its highest weight
vector. We always have $\sharp(\theta_k)\underset{M}\leq 
\flat(\theta_k)$.

\subsection{}
\label{(3)}
On line 5 in the fifth paragraph of Sect. 4.1 one must have 
$\sharp(\theta_k)$ instead of $\flat(\theta_k)$. 

\subsection{}
\label{(4)} In the proof of Lemma 4.3 one must everywhere replace
$\flat(\theta_k)$ by $\sharp(\theta_k)$.

\subsection{}
\label{(5)} In the proof of Lemma 4.3, the scheme 
$\underset{m}\Pi\, (\Gr_{M,X}^{\sharp(\theta_k)})^{\times n_m}$
must be replaced by the corresponding iterated convolution diagram.

In more detail, for $\theta=\underset{k}\Sigma\, \theta_k$, we regard
the $\{\theta_k\}$, $k=1,...,l$ (of course, $l=\underset{m}\Sigma\, n_m$)
as a set with multiplicities, and we choose an order on it.
We need the scheme that classifies the data of points
$(x_1,...,x_l)\in X$, principal $M$-bundles 
$\F^1_M,...,\F^{l+1}_M$ on $X$; and isomorphisms
$\beta_k:\F^k_M|_{X-x_k}\simeq \F^{k+1}_M|_{X-x_k}\in 
\ol{\H}_{M}^{\sharp(\theta_k)}$,
where $\ol{\H}_{M}^{\sharp(\theta_k)}$ is the corresponding
closed substack in the Hecke stack, i.e., the relative over
$\Bun_M$ version of $\ol{\Gr}_{M}^{\sharp(\theta_k)}$.

\subsection{}
\label{(6)} The claim in the first sentence in the proof of~Lemma~7.7 is wrong. Namely,
${\mathcal Q}_{{\mathfrak P}(\theta)}$ is not in general the tensor product of the
${\on{IC}}$-sheaf of the stratum $\overline{_{{\mathfrak P}(\theta)}{\on{\overline{Bun}}}_P}$ by a constant complex. Namely, the corresponding perverse sheaf is nonconstant along $X^{{\mathfrak P}(\theta)}$. This happens already for ${\operatorname{SL}}_n$ and a maximal parabolic attached to the root $e_1-e_2$, where your ${\on{\widetilde\Bun}}_P, {\on{\overline{Bun}}}_P$ are smooth.

\sssec{} Let $\Theta$ be the set of elements $\theta\in\Lambda_{G,P}$ such that $\check{{\mathfrak u}}(P)_{\theta}\ne 0$. Let ${\mathbb Z}_+\Theta$ be the free abelian group with base $\Theta$. The map $\Theta\to \Lambda_{G,P}$ extends to ${\mathbb Z}_+\Theta\to \Lambda_{G,P}$. For $\theta\in\Lambda_{G,P}$ and $V\in{\operatorname{Rep}}(\check{M})$ 
you denote by $V_{\theta}$ the direct summand on which $Z(\check{M})$ acts by $\theta$. 
An element of ${\mathbb Z}_+\Theta$ over $\theta\in\Lambda_{G,P}$ is denoted ${\mathfrak P}(\theta)$. 

\sssec{} Given ${\mathfrak P}(\theta)=\mathop{\sum}\limits_{\theta_s\in\Theta} n_s\theta_s\in{\mathbb Z}_+\Theta$ recall that $X^{{\mathfrak P}(\theta)}=\mathop{\prod}\limits_{\theta_s\in\Theta} X^{(n_s)}$. Say that a $Z(\check{M})$-graded module $V$ is $\Theta$-adapted if for any $\theta\in\Lambda_{G,P}$, $V_{\theta}=0$ unless $\theta\in\Theta$. 
For such ${\mathfrak P}(\theta)$ and a $\Theta$-adapted ${\mathbb Z}$-graded $Z(\check{M})$-module $V$ set for brevity
\[{\mathcal R}_{{\mathfrak P}(\theta)}(V)=\mathop{\boxtimes}\limits_{\theta_s\in \Theta} \, \overline{\Loc}(V_{\theta_s})^{(n_s)} \in Shv(X^{{\mathfrak P}(\theta)})\]
%Its restriction to the diagonal $X\hook{} X^{{\mathfrak P}(\theta)}$ is the constant complex 
%$$
%\bar{\mathcal R}_{{\mathfrak P}(\theta)}(V):=\mathop{\otimes}\limits_{\theta_s\in \Theta} \, \Sym^{n_s}\overline{\Loc}(V_{\theta_s})
%$$

For ${\mathfrak P}(\theta_1), {\mathfrak P}(\theta_2)\in {\mathbb Z}_+\Theta$ we write ${\mathfrak P}(\theta_1)+{\mathfrak P}(\theta_2)$ for their sum in the semigroup ${\mathbb Z}_+\Theta$. Let 
$$
s_{{\mathfrak P}(\theta_1), {\mathfrak P}(\theta_2)}: X^{{\mathfrak P}(\theta_1)}\times X^{{\mathfrak P}(\theta_2)}\to X^{{\mathfrak P}(\theta_1)+{\mathfrak P}(\theta_2)}
$$ 
be the  corresponding addition map. 

If now $V_1, V_2$ are $\Theta$-adapted ${\mathbb Z}$-graded $Z(\check{M})$-modules and $
{\mathfrak P}(\theta)\in{\mathbb Z}_+\Theta$ then
\begin{equation}
\label{iso_generality_about_gB}
{\mathcal R}_{{\mathfrak P}(\theta)}(V_1\oplus V_2)\,\stackrel{\sim}{\longrightarrow}\, \mathop{\oplus}\limits_{{\mathfrak P}(\theta_1), {\mathfrak P}(\theta_2)\in{\mathbb Z}_+\Theta, {\mathfrak P}(\theta_1)+{\mathfrak P}(\theta_2)={\mathfrak P}(\theta)} (s_{{\mathfrak P}(\theta_1), {\mathfrak P}(\theta_2)})_*({\mathcal R}_{{\mathfrak P}(\theta_1)}(V_1)\boxtimes {\mathcal R}_{{\mathfrak P}(\theta_2)}(V_2)
\end{equation}

\sssec{} For $\theta\in\Lambda_{G,P}^{pos}$ and ${\mathfrak P}(\theta)\in{\mathbb Z}_+\Theta$ over $\theta$ let $j_{P, {\mathfrak P}(\theta)}: X^{{\mathfrak P}(\theta)}\times{\on{\overline{Bun}}}_P\to {\on{\overline{Bun}}}_P$ denote the composition
\[X^{{\mathfrak P}(\theta)}\times{\on{\overline{Bun}}}_P\xrightarrow{i_{{\mathfrak P}(\theta)}\times\id} X^{\theta}\times{\on{\overline{Bun}}}_P\to 
{\on{\overline{Bun}}}_P\]
The map $j_{P, {\mathfrak P}(\theta)}$ is finite. Denote by $^{{\mathfrak P}(\theta)}{\on{\overline{Bun}}}_P$ its image. In fact, $j_{P, {\mathfrak P}(\theta)}$ is a normalization of $^{{\mathfrak P}(\theta)}{\on{\overline{Bun}}}_P$. So, 
\[(j_{P, {\mathfrak P}(\theta)})_!{\on{IC}}\,\stackrel{\sim}{\longrightarrow}\, {\on{IC}}(^{{\mathfrak P}(\theta)}{\on{\overline{Bun}}}_P).\] 

We have the principal nilpotent $e\colon \check{{\mathfrak u}}(P)\to \check{{\mathfrak u}}(P)$ defined
in~Section~7.1. Recall that ${\mathfrak r}: {\on{\widetilde\Bun}}_P\to{\on{\overline{Bun}}}_P$.
The complex ${\mathfrak Q}:={\mathfrak r}_!{\on{IC}}$ should be of the form
\begin{equation}
\label{complex_gQ_general}
{\mathfrak Q}\,\stackrel{\sim}{\longrightarrow}\, \mathop{\oplus}\limits_{\theta\in\Lambda_{G,P}^{pos},\;  
{\mathfrak P}(\theta)}
\; (j_{P, {\mathfrak P}(\theta)})_!({\mathcal R}_{{\mathfrak P}(\theta)}({\on{Im}} e)\boxtimes{\on{IC}}({\on{\overline{Bun}}}_P))[2|{\mathfrak P}(\theta)|]
\end{equation}

The proof is by induction going from the open stratum to more and more degenerate ones.

 Note that $({\on{Im}}\ e)[1]\in{\operatorname{Vect}}^c$ is self-dual. For this reason ${\mathcal R}_{{\mathfrak P}(\theta)}({\on{Im}}\ e)[2|{\mathfrak P}(\theta)|]$ on $X^{{\mathfrak P}(\theta)}$ is self-dual.

\bigskip

We would like to express our deep gratitude to S.~Lysenko
and D.~Nadler, who have pointed out mistakes corrected in 
subsections~\ref{(2)}--\ref{(6)} and~\ref{(1)} above, respectively.

\end{document}